\documentclass[11pt,leqno]{amsart}
\usepackage[T1]{fontenc}
\usepackage{ae}
\usepackage{aecompl}
\usepackage[cm]{aeguill}
\usepackage{amsmath}
\usepackage{amssymb}
\usepackage{latexsym,graphicx,chicago}
\usepackage{graphicx,multirow,rotating}

\newtheorem{thm}{Theorem}[section]
\newtheorem{lem}{Lemma}[section]
\newtheorem{prop}{Proposition}[section]
\newtheorem{cor}{Corollary}[section]

\newtheorem{rem}{Remark}[section]

\def\1{{{\mbox{${\rm{1\negthinspace\negthinspace I}}$}}}}
\newcommand{\eref}[1]{(\ref{#1})}

\newcommand\ind{{ {{1}}\hspace{-0,8mm}{\mathrm I}}}
\newcommand{\esp}{\mathbb{E}}

\newcounter{hypc}
\newcommand{\hypothese}[2]{\stepcounter{hypc}
\tag{$\mathbf{A}_{\thehypc}^{#2}$}
\label{#1}}
\newcounter{cond}

\begin{document}
\title[Adaptive density deconvolution with dependent inputs]{Adaptive density deconvolution with dependent inputs}
\author{F. Comte$^{*,1}$}\thanks{$^1$ Universit\'e Paris V,  MAP5,
UMR CNRS 8145.}
\address{F. Comte\\ MAP5 UMR 8145 \newline
Universit\'e Ren\'e Descartes-Paris 5\\ 45 rue des Saints-P\`eres
75270 Paris cedex 06, France.\newline email:
fabienne.comte@univ-paris5.fr.}
\author{J. Dedecker$^2$}\thanks{$^2$ Universit\'e Paris 6, Laboratoire de
Statistique Th\'eorique et Appliqu\'ee.}
\address{J. Dedecker\\
 Laboratoire de Statistique Th\'eorique et Appliqu\'ee \newline
Universit\'e paris 6, 175 rue du Chevaleret  75013 Paris, France.
\newline email:
dedecker@ccr.jussieu.fr}
\author{M. L. Taupin $^3$}\thanks{$^3$ IUT de Paris V et
Universit\'e Paris Sud, Laboratoire de Probabilit\'es, Statistique
et Mod\'elisation, UMR 8628.}
\address{M.-L. Taupin\\
Laboratoire de Probabilit\'es, Statistique et
Mod\'elisation, UMR C 8628 \newline Universit\'e Paris-Sud\\
B\^{a}timent 425\\ 91405 Orsay Cedex, France.
\newline email:
Marie-Luce.Taupin@math.u-psud.fr}
\begin{abstract}
In the convolution model $Z_i=X_i+ \varepsilon_i$, we give a model
selection procedure to estimate the density of the unobserved
variables $(X_i)_{1 \leq i \leq n}$, when the sequence $(X_i)_{i
\geq 1}$ is strictly stationary but not necessarily independent.
This procedure depends on wether the density of $\varepsilon_i$ is
super smooth or ordinary smooth. The rates of convergence of the
penalized contrast estimators are the same as in the independent
framework, and are minimax over most classes of regularity on
${\mathbb R}$. Our results apply to mixing sequences, but also to
many other dependent sequences. When the errors are super smooth,
the condition on the dependence coefficients is the minimal
condition of that type ensuring that the sequence $(X_i)_{i \geq 1}$
is not a long-memory process.
\end{abstract}
\maketitle
\begin{center}
\end{center}

Deconvolution. Dependence. Mixing. Penalized Contrast. Hidden Markov models

\bibliographystyle{chicago}
\section{Introduction}
\setcounter{equation}{0} \setcounter{lem}{0}
\setcounter{thm}{0}

The problem of estimating the density of identically distributed but
not independent random variables $X_1, \dots, X_n$ when they are
observed with an additive and independet noise is encountered in numerous contexts. This
problem is described by the model
\begin{eqnarray}
Z_i=X_i+\varepsilon_i, \quad \text{ for } i=1, \dots,
n,\label{model}\end{eqnarray} where one observes $Z_1,\ldots,Z_n$,
and where  $(\varepsilon_i)_{1 \leq i \leq n}$ are independent and
identically distributed (i.i.d.), and independent of $(X_i)_{1 \leq
i \leq n}$. When $(X_i)_{i \leq 1 \leq n}$ is a Markov
chain, the model \eref{model} is a particular case of hidden Markov
models, with an additive structure.

 Our aim is the adaptive estimation of $g$, the
common  distribution of the unobserved variables $(X_i)_{1 \leq i
\leq n}$, when the density $f_\varepsilon$ of $\varepsilon_i$ is
known. More precisely we shall build an estimator  of $g$ without
any prior knowledge on its smoothness, using the observations
$(Z_i)_{i \leq 1 \leq n}$ and the knowledge of the convolution
kernel $ f_{\varepsilon}$.
We shall assume that the known density $f_{\varepsilon}$  belongs to
various collections of densities, and that the dependence properties
of the sequence $(X_i)_{i \geq 1}$ are described by appropriate
dependence coefficients. More precisely, we consider two types of
dependent sequences. We assume either that the sequence $(X_i)_{i
\geq 1}$ is absolutely regular in the sense of Rozanov and
Volkonskii \citeyear{VolkRozan}, or that it is   $\tau$-dependent in
the sense of Dedecker and Prieur~\citeyear{JDCPCoeff}. These
dependence  conditions are presented in Section \ref{dep}  and
motivated through various examples.

In density deconvolution, two factors determine the estimation
accuracy. First, the smoothness of the density $g$ to be estimated,
and second the smoothness of the error density, the worst rates of
convergence being obtained for the smoothest errors densities.
We shall consider two classes of densities for $f_\varepsilon$:
first the so called super smooth densities  with exponential  decay
of their Fourier transform, and next the class  of ordinary smooth
densities with Fourier transform having a polynomial decay.

Let us briefly recall the  previous results in the independent
framework. To our knowledge, the first adaptive estimator has been
proposed by Pensky and Vidakovic~\citeyear{Penskyvidakovic99}. It is
a wavelet estimator constructed $via$ a thresholding procedure. This
estimator achieves the minimax rates when $g$ belongs to a
Sobolev class, but it fails to reach the minimax rates when both the
errors density and $g$ are supersmooth.
More recently, Comte {\it et al.}~\citeyear{CRT1} have proposed an
adaptive estimator of $g$ constructed by   minimizing an appropriate
penalized contrast function only depending on the observations and
on $f_\varepsilon$. This estimator is minimax (sometimes within a
negligible logarithmic factor) in all cases where lower bounds are
previously known (i.e.  in most cases). More precisely, the authors
obtain non-asymptotic upper bounds for the Mean Integrated Square
Error (MISE), which ensure an automatic trade-off between a bias
term and the penalty term. Hence, the estimator automatically
achieves the best rate obtained by the collection of non-penalized
estimators when the (unknown) optimal space is selected (sometimes
up to a negligible logarithmic factor).
When both the density and the errors are super smooth, this adaptive
estimator significantly improves on the rates given by the adaptive
estimator built in Pensky and
Vidakovic~\citeyear{Penskyvidakovic99}, whereas both adaptive
estimators have the same rate in the other cases. This improvement
partly comes from the choice of the Shannon basis (see Section
\ref{Shannon}) instead of the wavelet basis considered in Pensky and
Vidakovic.


In the dependent context, we follow the approach  proposed in Comte
{\it et al.} \citeyear{CRT1}.  We give adaptive estimators of $g$,
constructed by  minimizing  an appropriate penalized contrast
function. The penalty function depends on the known density
$f_\varepsilon$, but it  does not depend on the dependence
coefficients of the sequence $(X_i)_{i \geq 1}$. The adaptive
estimators have the same rates as in the independent case, under
mild conditions on the dependence coefficients of $(X_i)_{i \geq
1}$. The important point here is that the penalty functions are the
same (or almost the same) as in the independent framework. This is a
bit surprising: indeed, when the  $(X_i)_{1 \leq i \leq n}$ are
observed (i.e. $\varepsilon_i=0$), the threshold level proposed in
Tribouley and Viennet~\citeyear{TribouleyViennet} as well as the
penalty function given in Comte and
Merlev\`ede~\citeyear{ComteMerlevede} (see also our Corollary
\ref{densitycor}) depends on the mixing coefficients of the sequence
$(X_i)_{i \geq 1}$.

In Section \ref{sec4} we deal with non adaptive estimators. As
usual, we show that the MISE of the minimum contrast estimator is
bounded by a squared bias plus a variance term. The variance term
can be split into two terms. The first and dominating term of the
variance is exactly  the variance of a density deconvolution
estimator in the independent context. It is as usual related to
$\int_{\vert  x\vert \leq C_n} \vert f_\varepsilon^*(x)\vert^{-2}
dx$, $C_n\rightarrow \infty$. The second  and negligible  term in
the variance is the term  involving the dependence structure of the
sequence $(X_i)_{i \geq 1}$.  The main consequence of this first
result is that this non adaptive estimator reaches the (minimax)
rates of the i.i.d. case (as given in Fan~\citeyear{fan91},
Butucea~\citeyear{Butucea}, and  Butucea and
Tsybakov~\citeyear{ButTsyb}), as soon as the dependence coefficients
are summable. Moreover, even if the coefficients are not summable,
there is no loss in the rate provided that the partial sums of the
coefficients does not grow too fast with respect to $\int_{\vert
x\vert \leq C_n} \vert f_\varepsilon^*(x)\vert^{-2} dx$. These results have to
be  compared with  previously known results for non adaptive density deconvolution in
dependent contexts. For strongly mixing sequences in the sense of Rosenblatt
\citeyear{Rosenblatt}, Masry \citeyear{Masry932} propose a
kernel-type estimator for the joint density $g_p$ of $(X_1, \ldots ,
X_p)$ when it exists. For the (pointwise) Mean Square Error, he
obtains the same rates as in the i.i.d. case provided that $\alpha(n)=
O(n^{-2-\delta})$ for ordinary smooth $f_\varepsilon$, and provided
that $\alpha(n)= O(n^{-1-\delta})$ for super smooth $f_\varepsilon$.
When $p=1$, our assumption on the mixing coefficients is
weaker, since we only need $\sum_{n>0} \alpha(n) < \infty$ in both
cases (see our Remark \ref{strongmix}).


In the  main part (Section \ref{sec5}), we study the adaptive
estimators. We show that the squared bias term and the variance term
obtained in the upper bound of the MISE of the adaptive estimator
are the same as in the independent case. The model selection
procedure depends on wether the density $f_\varepsilon$ is super
smooth or ordinary smooth.

When  $f_\varepsilon$ is super smooth, the adaptive estimator, is
constructed with the exact penalty of the independent context. Its
rate of convergence is exactly the same as in the independent case,
provided that the dependence coefficients of $(X_i)_{i \geq 1}$ are
summable. The main tools in this case are covariance inequalities
for dependent variables, and concentration inequalities. The case of
super smooth errors is particularly important, since it
contains the case of Gaussian errors. It also contains the
stochastic volatility model, in which $\varepsilon_i\sim
\ln({\mathcal N}(0,1)^2)$ (see  Van Es {\it et
al.}~\citeyear{vanes03,vanes05}, Comte~\citeyear{comtejtsa}, Comte and
Genon-Catalot \citeyear{CGC}).

When $f_\varepsilon$ is ordinary smooth, the adaptive estimator, is
constructed with a penalty of the same order as in the independent
context.
 Its rate of convergence is exactly the same as in the independent case.
For ordinary smooth errors, the main tools are the coupling
properties of the  dependence coefficients (see Section
\ref{couplage}). To use these  properties, we need to consider a
more restrictive type of dependence than for super smooth errors,
and we need to impose a polynomial decrease of the coefficients.

In both cases, super and ordinary smooth, the results hold for
$\beta$-mixing and $\tau$-dependent  random variables $(X_i)_{i \geq
1}$. To our knowledge, this is the first time that adaptive density
deconvolution in a dependent context is considered. The robustness
of this estimation procedure to dependency strongly use the
independence between $(X_i)_{1 \leq i \leq n} $ and
$(\varepsilon_i)_{i \leq 1 \leq n}$, and the fact that the errors
are i.i.d. random variables.
We refer to Comte {\it et al.}~\citeyear{CRT2,CRT1} for practical
implementation of the estimators, and  for the calibration of the
constants in the penalty functions. In Comte {\it et
al.}~\citeyear{CRT2}, the robustness of the procedure to various
dependency has been experimented in practice (see Tables 4 and 5
therein).


\section{Some measures of dependence} \label{dep}
\setcounter{equation}{0}
\setcounter{lem}{0}
\setcounter{thm}{0}
Let $(\Omega, {\mathcal A}, {\mathbb P})$ be a probability space.
Let $Y$ be a random variable with values in a Banach space
$({\mathbb B}, \|\cdot \|_{\mathbb B})$, and let ${\mathcal M}$ be a
$\sigma$-algebra of ${\mathcal A}$. Let ${\mathbb P}_{Y|{\mathcal
M}}$ be a conditional distribution of $Y$ given ${\mathcal M}$, and
let $P_Y$ be the distribution of $Y$. Let ${\mathcal B}({\mathbb
B})$ be the borel $\sigma$-algebra on $({\mathbb B}, \| \cdot
\|_{\mathbb B})$, and let $\Lambda_1({\mathbb B})$ be the set of
1-Lipschitz functions from $({\mathbb B}, \|\cdot\|_{\mathbb B})$ to
${\mathbb R}$. Define now
\begin{eqnarray*}
\beta({\mathcal M}, \sigma(Y))&=& {\mathbb E}\Big (\sup_{A \in
{\mathcal B}({\mathcal X})} |{\mathbb P}_{Y|{\mathcal
M}}(A)-{\mathbb P}_Y(A)|
\Big ) \, , \\
\text{and if ${\mathbb E}(\|Y\|) < \infty$,} \quad \tau ({\mathcal
M}, Y)&=& {\mathbb E}\Big (\sup_{f \in {\Lambda_1}({\mathbb B})}
|{\mathbb P}_{Y|{\mathcal M}}(f)-{\mathbb P}_Y(f)| \Big ) \, .
\end{eqnarray*}
The coefficient $\beta({\mathcal M}, \sigma(Y))$ is the usual
mixing coefficient, introduced by Rozanov and Volkonskii
\citeyear{VolkRozan}. The coefficient $\tau({\mathcal M}, Y)$ has
been introduced by Dedecker and Prieur \citeyear{JDCPCoeff}.

Let ${\bf X}=(X_i)_{i \geq 1}$ be a strictly stationary sequence
of real-valued random variables. For any $k\geq 0$, the
coefficients $\beta_{{\bf X},1}(k)$ and $\tau_{{\bf X}, 1}(k)$ are
defined by
\begin{eqnarray*}
\beta_{{\bf X}, 1}(k)&=&  \beta(\sigma (X_1), \sigma(X_{1+k})), \\
\text{and if ${\mathbb E}(|X_1|) < \infty$,} \quad \tau_{{\bf
X},1}(k)&=&  \tau(\sigma (X_1), X_{1+k}).
\end{eqnarray*}
On ${\mathbb R}^l$, we put the norm $\|x-y\|_{{\mathbb
R}^l}=l^{-1}(|x_1-y_1|+ \cdots +|x_l-y_l|)$. Let ${\mathcal
M}_i=\sigma (X_k, 1 \leq k \leq i)$. The coefficients $\beta_{{\bf
X},\infty}(k)$ and $\tau_{{\bf X}, \infty}(k)$ are defined by
\begin{eqnarray*}
\beta_{{\bf X}, \infty}(k)&=& \sup_{i \geq 1, l \geq 1} \sup \left
\{ \beta({\mathcal M}_i, \sigma(X_{i_1}, \ldots , X_{i_l})), i+k
\leq i_1 < \cdots < i_l \right \},
\\
\text{and if ${\mathbb E}(|X_1|) < \infty$,} \quad \tau_{{\bf X},
\infty}(k)&=& \sup_{i \geq 1, l \geq 1} \sup \left \{
 \tau({\mathcal M}_i, (X_{i_1}, \ldots , X_{i_l})), i+k \leq i_1 < \cdots < i_l \right \}. \\
\end{eqnarray*}
Let $Q_X$ be the generalized inverse of the tail function $x
\rightarrow {\mathbb P}(|X_1|>x)$. We have the inequalities
\begin{eqnarray}
\tau_{{\bf X}, 1}(k) \leq 2 \int_0^{\beta_{{\bf X}, 1}(k)} Q_X(u)
du  &\text{and} & \quad \tau_{{\bf X}, \infty}(k) \leq 2 \int_0^{
\beta_{{\bf X}, \infty}(k)} Q_X(u) du \, .
\end{eqnarray}

\subsection {Coupling}
\label{couplage}
We recall the coupling properties of these
coefficients. Assume that $\Omega$ is rich enough, which means that
there exists $U$ uniformly distributed over $[0,1]$ and independent
of ${\mathcal M} \vee \sigma (X)$. There exist two ${\mathcal M}
\vee \sigma(U) \vee \sigma (X)$-measurable random variables $X_1^*$
and $X_2^*$ distributed as $X$ and independent of ${\mathcal M}$
such that
\begin{equation}\label{cber}
\beta ({\mathcal M}, \sigma(X))= {\mathbb P}(X \neq X_1^*) \quad
\text{and} \quad  \tau({\mathcal M}, X)= {\mathbb E}(\|X-
X_2^*\|_{\mathbb B}) \, .
\end{equation}
The first equality in (\ref{cber}) is due to Berbee \citeyear{Berbee}, and the
second one has been established in Dedecker and Prieur \citeyear{JDCPCoeff},
Section 7.1.

\subsection {Covariance inequalities}
Denote by $\| \cdot \|_{\infty, {\mathbb P}}$ the ${\mathbb
L}^\infty (\Omega, {\mathbb P})$-norm. Let $X, Y$ be two real-valued
random variables,  and let $f, h$ be two measurable functions from
${\mathbb R}$ to ${\mathbb C}$. Then
\begin{equation}\label{ibrb}
|\text{Cov}(f(Y), h(X))| \leq 2 \|f(Y)\|_{\infty, {\mathbb P}}
\|h(X)\|_{\infty, {\mathbb P}} \,  \beta (\sigma (X), \sigma (Y)) \,
,
\end{equation}
and if $\text {Lip}(h)$ is the Lipschitz coefficient of $h$,
\begin{equation}\label{evident}
 \quad |\text{Cov }(f(Y), h(X))|
\leq \|f(Y)\|_{\infty, {\mathbb P}} \text{Lip}(h) \, \tau (\sigma
(Y), X) \, .
\end{equation}
Inequalities \eref{ibrb} and \eref{evident} follow from the
coupling properties \eref{cber} by noting that if $X^*$ is
distributed as $X$ and independent of $Y$,
$$ \text{Cov }(f(Y), h(X))= {\mathbb E}(\overline{f(Y)}(h(X)-h(X^*)))\, .$$

\subsection {Examples} Examples of $\beta$-mixing sequences are well known (we refer to the books by
Doukhan \citeyear{Doukhan} and Bradley \citeyear{Bradley}).
One of the most important
examples is
the following: a stationary, irreducible,
aperiodic and positively recurrent Markov chain $(X_i)_{i \geq 1}$
is $\beta$-mixing, which means that $\beta_{{\bf X}, \infty}(k)$
tends to zero as $k$ tends to infinity.

Unfortunately, many simple Markov chains are not $\beta$-mixing (and
not even  strongly mixing in the sense of Rosenblatt
\citeyear{Rosenblatt}). For instance, if $(\epsilon_i)_{i \geq 1}$
is i.i.d. with marginal ${\mathcal B}(1/2)$, then the stationary
solution $(X_i)_{ i \geq 0}$ of the equation
\begin{equation}
                  X_n= \frac{1}{2} (X_{n-1}+ \epsilon_n),  \quad
                  \text{$X_0$ independent of $(\epsilon_i)_{i \geq
                  1}$}
\end{equation}
is not $\beta$-mixing (and not even  strongly mixing) since
$\beta_{{\bf X},1}(k)=1$ for any $k\geq 0$. By contrast, for this
particular example, one has $\tau_{{\bf X}, \infty}(k) \leq 2^{-k}$.
More generally, the coefficient $\tau_{{\bf X}, \infty}(k)$ is easy
to compute in many situations (see Dedecker and Prieur
\citeyear{JDCPCoeff}). Let us recall some important examples:

\medskip

\noindent {\bf Linear processes.} Assume that $X_i=\sum_{j \geq 0}
a_j \xi_{n-j}$, where $(\xi_i)_{i \in {\mathbb Z}}$ is i.i.d. One
has the bounds
$$
  \tau_{{\bf X}, \infty}(k) \leq 2 {\mathbb E}(|\xi_0|) \sum_{j \geq k} |a_j|
  \quad \text{and} \quad
  \tau_{{\bf X}, \infty}(k) \leq \sqrt{ 2 \text{Var}(\xi_0) \sum_{j
  \geq k} a_j^2} .
$$

\noindent {\bf Markov chains.} Let $(X_n)_{n \geq 0}$ be a
stationary Markov chain such that $X_n=F(X_{n-1}, \xi_n)$ for some
measurable function $F$ and some i.i.d. sequence $(\xi_i)_{i \geq 1}$
independent of $X_0$. Assume that there exists $\kappa <1$ such that
$$
 {\mathbb E}( |F(x, \xi_0)-F(y, \xi_0)|) \leq \kappa |x-y| \, .
$$
Then one has the inequality
$$
\tau_{{\bf X}, \infty}(k) \leq 2  {\mathbb E}(|X_0 |) \kappa^k\, .
$$
An important example is $X_n=f(X_{n-1}) + \xi_n$ for some
$\kappa$-lipschitz function $f$.

\medskip

\noindent{\bf Expanding maps.} Let $T$ be a Borel-measurable map
from $[0,1]$ to $[0,1]$. If the probability $\mu$ is invariant by
$T$, the sequence $(Y_i=T^i)_{i \geq 0}$ of random variables from
$([0,1], \mu)$ to $[0,1]$ is strictly stationary. Define the
operator $K$ from ${\mathbb L}^1([0,1], \mu)$ to ${\mathbb
L}^1([0,1], \mu)$ $via$ the equality
$$
\int_0^1 (Kh)(x) k(x)  \mu (dx) = \int_0^1h(x) (k \circ T)(x) \mu
(dx)
$$
where $h \in {\mathbb L}^1([0,1], \mu)$ and $k \in {\mathbb
L}^\infty([0,1], \mu)$. It is easy to check  that $(Y_1, Y_2, \ldots
, Y_n)$ has the same distribution as $(X_n,X_{n-1},\ldots,X_1)$
where $(X_i)_{i \in {\mathbb Z}}$ is a stationary Markov chain with
invariant distribution $\mu$ and transition kernel $K$. If $T$ is
uniformly expanding (see for instance the assumptions on page 218 in
Dedecker and Prieur \citeyear{JDCPCoeff}), then there exist $C>0$
and $\rho $ in $]0,1[$ such that
$$
   \tau_{\bf {X}, \infty} (k) \leq C \rho^k
$$
(see Dedecker and Prieur page 230). Note that the Markov chain
$(X_i)_{i \geq 1}$ is not $\beta$-mixing (and not even strongly
mixing). Indeed $\beta( \sigma(X_1), \sigma (X_n))= \beta (\sigma
(T^n), \sigma (T))$. Since $\sigma (T^n) \subset \sigma (T)$, it
follows that
$$\beta( \sigma(X_1), \sigma (X_n))\geq \beta (\sigma (T^n), \sigma
(T^n))= \beta (\sigma(T), \sigma (T))$$ and the later is positive as
soon as $\mu$ is non trivial.

\section{Assumptions and estimators}
\setcounter{equation}{0}
\setcounter{lem}{0}
\setcounter{thm}{0}
\bigskip\noindent
For two complex-valued functions $u$ and $v$ in
$\mathbb{L}_2(\mathbb{R})\cap \mathbb{L}_1(\mathbb{R})$, let
$$
u^*(x)=\int e^{itx}u(t)dt,  \quad u*v(x)=\int u(y)v(x-y)dy, \
\text{and} \quad  <u,v>=\int u(x)\overline{v}(x)dx$$ with
$\overline{z}$ the conjugate of a complex number $z$. We also use
the notations $$ \|u\|_1=\int |u(x)|dx, \quad \|u\|^2=\int |u(x)|^2
dx, \quad \text{and} \quad \|u\|_\infty=\sup_{x \in
\mathbb{R}}|u(x)|.
$$

\subsection{Assumptions for density deconvolution}\label{add}
The smoothness of $f_\varepsilon$ is described by the following assumption.
\begin{align}
& \mbox{ There exist nonnegative numbers } \kappa_0, \gamma,~ \mu, \mbox{ and }\delta \mbox{ such
that $f_{\varepsilon}^*$ satisfies }\notag\\
& \kappa_0(x^2+1)^{-\gamma/ 2}\exp\{-\mu\vert x\vert^\delta\}\leq
|f_\varepsilon^*(x)| \leq 
\kappa_0'(x^2+1)^{-\gamma/ 2}\exp\{-\mu\vert x\vert^\delta\}.
\hypothese{condfeps}{\varepsilon}\\
&\hypothese{fepsnn}{\varepsilon}\mbox{ The density }f_\varepsilon \mbox{
belongs to } \mathbb{L}_2(\mathbb{R})
\mbox{ and for all }x\in
\mathbb{R},\,f_\varepsilon^*(x)\not=0.
\end{align}

Since $f_\varepsilon$ is known, the constants $\mu, \delta,
\kappa_0,$ and $\gamma$  defined in \eref{condfeps}  are also known.

When $\delta=0$ in \eref{condfeps}, $f_\varepsilon$ is  usually
called ``ordinary smooth''. When $\mu>0$ and $\delta>0$,
$f_\varepsilon$ is called ``super smooth''. Densities satisfying
\eref{condfeps} with $\delta>0$ and $\mu>0$ are infinitely
differentiable. The standard examples for super smooth densities are
the following: Gaussian or Cauchy distributions are super smooth of
order $\gamma=0, \delta=2$ and $\gamma=0, \delta=1$ respectively.
When $\varepsilon=\ln(\eta^2)$ with $\eta\sim {\mathcal N}(0,1)$ as
in Van Es {\it et al.} \citeyear{vanes03,vanes05}, then
$\varepsilon$ is super-smooth with $\delta=1, \gamma=0$ and $\mu=
\pi/2$. For ordinary smooth densities, one can cite for instance the
double exponential (also called Laplace) distribution with
$\delta=0=\mu$ and $\gamma=2$. Although densities with $\delta>2$
exist, they are difficult to express in a closed form. Nevertheless,
our results hold for such densities. Furthermore, the square
integrability of $f_{\varepsilon}$ in \eref{fepsnn}  require that
$\gamma> 1/2$ when $\delta=0$ in \eref{condfeps}.

Classically,  the slowest rates of convergence for estimating $g$
are obtained for super smooth error densities. In particular, when
$\varepsilon$ is Gaussian and $g$ belongs to Sobolev classes, the
minimax rates are negative powers of $\ln(n)$ (see
Fan~\citeyear{fan91}). Nevertheless, the rates are improved if $g$
has stronger smoothness properties, described  by the set
\begin{eqnarray}
\label{super} \mathcal{S}_{s,r,b}(C_1)=\Big \{\psi\; \mbox{ such that }
\int_{-\infty}^{+\infty} |\psi^*(x)|^2(x^2+1)^{s}\exp\{2b |x|^{r} \}
dx\leq C_1 \Big \}
\end{eqnarray}
for $s,r,b$ non-negative numbers.

Such smoothness classes are classically considered both in
deconvolution and in density estimation without errors. When $r=0$,
\eref{super} corresponds to a Sobolev ball.
 The functions in \eref{super} with $r>0$ and  $b>0 $ are infinitely many times differentiable.
They admit analytic continuation on a finite width strip when $r=1$
and on the whole complex plane if $r=2$.

Subsequently, the density $g$ is supposed to satisfy the following
assumption.
\begin{align}
& \mbox{The density } g\in \mathbb{L}_2(\mathbb{R}) \mbox{
  and there exists }  M_2>0, \mbox{ such that }
\int x^2g^2(x)dx <M_2<\infty \hypothese{momg}{X}.
\end{align}
Assumption \eref{momg} which is due to the construction of the
estimator,  is quite unusual in density estimation. It already
appears in density deconvolution in the  independent framework in
Comte {\it et al.}~\citeyear{CRT2,CRT1}. It also appears in a
slightly different way in Pensky and
Vidakovic~\citeyear{Penskyvidakovic99} who assume, instead of
\eref{momg}  that $\sup_{x \in \mathbb{R}} \vert x\vert
g(x)<\infty$. It is important to note that Assumption \eref{momg} is
very unrestrictive.

All densities having tails of order $|x|^{-(s+1)}$ as $x$ tends to
infinity satisfy \eref{momg}  only if $s>1/2$. One can cite for
instance the Cauchy distribution or all stable distributions with
exponent $r>1/2$ (see Devroye~\citeyear{Devroye86}). The L\'evy
distribution, with exponent $r=1/2$ does not satisfies \eref{momg}.

\subsection{The projection spaces } \label{Shannon}

Let $\varphi(x)=\sin(\pi x)/(\pi x)$. For $m\in {\mathbb N}$ and
$j\in {\mathbb Z}$, set $\varphi_{m,j}(x) =\sqrt{m}
\varphi(mx-j).$ The functions $\{\varphi_{m,j}\}_{j \in
\mathbb{Z}}$ constitute an orthonormal system in ${\mathbb
L}^2({\mathbb R})$ (see e.g.  Meyer \citeyear{MeyerI}, p.22).  For
$m=2^k$, it is known as the Shannon basis. Though we choose here
integer values for $m$, a thinner grid would also be possible. Let
us define
$$S_m= \overline{{\rm span}}\{\varphi_{_{m,j}}, \;
j\in \mathbb{Z}\}, \ m\in {\mathbb N}.$$ The space $S_m$ is exactly
the  subspace of ${\mathbb L}_2({\mathbb R})$ of functions having a
Fourier transform with compact support contained in $[-\pi m, \pi
m]$.

The orthogonal projections of $g$ on $S_m$ is $g_m=\sum_{j\in
{\mathbb Z}} a_{m,j}(g) \varphi_{m,j}$ where $a_{m,j}(g) =
<\varphi_{m,j},g>$. To obtain representations having a finite number
of  ``coordinates'', we introduce
$$S_{m}^{(n)}= \overline{\rm span }\left\{\varphi_{m,j}, |j|\leq
k_n\right\}$$ with integers
$k_n$ to be specified later. The family
$\{\varphi_{m,j}\}_{\vert j\vert \leq k_n}$ is an orthonormal
basis of $S_m^{(n)}$ and the orthogonal projections of $g$ on $S_{m}^{(n)}$ is
 given by $g_{m}^{(n)}=\sum_{|j|\leq k_n}
a_{m,j}(g) \varphi_{m,j}$.

\subsection{Construction of the minimum contrast estimators}
\label{nonpen}

\bigskip\noindent For an arbitrary fixed integer $m$, an estimator of $g$
belonging to $S_m^{(n)}$ is defined by
\begin{equation}\label{tronque}  \hat g_m^{(n)} = \arg\min_{t\in S_m^{(n)}}
\gamma_n(t),\end{equation} where, for $t$ in $S_m^{(n)}$,
\begin{equation*}\label{u*}
\gamma_n(t)=\frac{1}{n}\sum_{i=1}^n \left[\|t\|^2
-2u_t^*(Z_i)\right], \;\;\; \mbox{ with } \;\;\;u_t(x) = \frac 1{2
\pi} \left(\frac{t^*(-x)}{f_{\varepsilon}^*(x)}\right).
\end{equation*}
By using Parseval and inverse Fourier formulae we obtain  that
$\mathbb{E}\left[u_t^*(Z_i)\right] =\langle t, g\rangle,$ so that
$\mathbb{E}(\gamma_n(t))=\|t-g\|^2 -\|g\|^2 $ is minimal when $t=g$.
This shows that $\gamma_n(t)$ suits well for the estimation of $g$.
Classical calculations show that $$\hat g_m^{(n)} = \sum_{\vert j\vert\leq
k_n} \hat a_{m,j} \varphi_{m,j} \;\; \mbox{with}\;\; \hat a_{m,j}=
\frac{1}{n} \sum_{i=1}^n u_{\varphi_{m,j}}^*(Z_i), \quad \text{and}
\quad \mathbb{E}(\hat a_{m,j})= <g,\varphi_{m,j}>=a_{m,j}.$$

\subsection{Minimum penalized contrast estimator}

As in the independent framework, the minimum penalized estimator
of $g$ is defined  as $\tilde g=\hat g_{\hat m_g}$ where $\hat
m_g$ is chosen in a purely data-driven way. The main point of the
estimation procedure lies in the choice of $m=\hat m_g$ for the
estimators $\hat g_m$ from Section \ref{nonpen} in order to mimic
the oracle parameter
\begin{eqnarray}
\label{oraclepar}
\breve m_g=\arg\min_{m}\mathbb{E}\parallel \hat  g_m-g\parallel_2^2.
\end{eqnarray}
The model selection is performed in an automatic way, using the
following penalized criteria
\begin{equation}\label{estitronc}
\tilde g=\hat g^{(n)}_{\hat m} \mbox{ with } \hat m= \arg\min_{m\in
\{1, \cdots , m_n \}} \left[\gamma_n(\hat g_m^{(n)}) + \; {\rm
pen}(m)\right],
\end{equation}
where $\mbox{pen}(m)$ is a penalty function, precised in the
Theorems, that depends on $f_\varepsilon^*$ through $\Delta(m)$
defined by
\begin{eqnarray}
\Delta(m)=\frac{1}{2\pi}\int_{-\pi m}^{\pi m}\frac{1}{\vert
f_\varepsilon^*( x)\vert^2}dx.\label{Delta1}
\end{eqnarray}
The key point in the dependent context is to find a penalty
function not depending on the mixing coefficients such that
$$\mathbb{E}\parallel \tilde g -g\parallel^2\leq C\inf_{m\in
\{1,\cdots,m_n\}}\mathbb{E}\parallel \hat g_m-g\parallel^2.$$


%
%
%
%
%
%
%
%
%


\section{Risk bounds for the minimum contrast estimators $\hat
g_m^{(n)}$} \label{sec4}
\setcounter{equation}{0} \setcounter{lem}{0}
\setcounter{thm}{0} We focus here on non adaptive estimation,
starting with the presentation of general upper bounds for MISEs of
the minimum contrast estimators $\hat g_m^{(n)}$.

\begin{prop}
\label{Vitss} If  \eref{fepsnn} and \eref{momg} hold, then
$$ {\mathbb E}\|g-\hat g_m^{(n)}\|^2\leq
\|g-g_m\|^2+ \frac{m^2(M_2+1)}{k_n} + \frac{2\Delta(m)}{n}+
\frac{2R_{m}}{n},
$$
where \begin{equation}\label{Rm}
 R_{m}=\frac{1}{\pi}\sum_{k=2}^{n}\int_{-\pi
m}^{\pi m}\left\vert {\mathrm
{Cov}}\left(e^{ixX_1},e^{ixX_k}\right)\right\vert dx.
\end{equation}
Moreover, $R_m \leq \min ( R_{m, \beta},  R_{m,
\tau})$, where
\begin{equation*}
  R_{m, \beta}= 4 m \sum_{k=1}^{n-1}\beta_{\mathbf{X},1}(k)
\quad \text{and} \quad
  R_{m, \tau} = \pi m^2 \sum_{k=1}^{n-1}\tau_{\mathbf{X},1}(k) \, .
\end{equation*}

\begin{rem}\label{strongmix}
{\rm The term $R_m$ can be easily bounded for many other dependent
sequences. For instance, if $\alpha_{\mathbf{X}, 1}= \alpha (\sigma
(X_1), \sigma (X_{1+k}))$ is the usual strong mixing coefficient of
Rosenblatt  \citeyear{Rosenblatt}, one has the upper bound $R_m \leq
16 m \sum_{k=1}^{n-1}\alpha_{\mathbf{X},1}(k)$. If ${\mathbf{X}}$ is
a stationary  sequence of associated random variables (see Esary
{\it et al.} \citeyear{Esaryetal} for the definition), then
$|{\mathrm {Cov}}(e^{ixX_1},e^{ixX_k})| \leq 4 x^2 {\mathrm
{Cov}}(X_1, X_k)$, so that $R_m \leq (8\pi^2/3) m^3
\sum_{k=2}^{n}{\mathrm {Cov}}(X_1, X_k)$.
 For more about density deconvolution with associated inputs,
we refer to the paper by Masry \cite{Masry03}.}
\end{rem}

\end{prop}
We now comment the rates resulting from Proposition \ref{Vitss}. As usual, the variance term $n^{-1}\Delta(m)$
depends on the rate of decay of the Fourier transform of
$f_\varepsilon$. According to Lemma \ref{details} and according to
Butucea and Tsybakov \citeyear{ButTsyb}, under
\eref{condfeps}-\eref{fepsnn}, we have
\begin{eqnarray}\nonumber
\lambda_1(f_\varepsilon,\kappa_0')\Gamma(m)(1+o(1))\!\!\!&\leq&\!\!\!
\Delta(m)\leq \lambda_1(f_\varepsilon,\kappa_0)\Gamma(m)(1+o(1))
\quad \text{as $m\rightarrow\infty$} \\\label{gammadem} \text{ where
} \quad \Gamma(m)\!\!\!&=&\!\!\! (1+(\pi m)^2)^{\gamma}(\pi
m)^{1-\delta} \exp\left\{2\mu(\pi m)^{\delta}\right\},
\end{eqnarray}
\begin{eqnarray}
\label{lambda1}
\qquad\quad\lambda_1(f_\varepsilon,\kappa_0)=\frac{1}
{\kappa_0^{2}\pi R(\mu,\delta)}, \mbox{ and }\, R(\mu,\delta)=
\mbox{1}\!\mbox{I}_{\{\delta=0\}} + 2\mu\delta
\mbox{1}\!\mbox{I}_{\{\delta>0\}} .
\end{eqnarray}
If  \eref{condfeps}-\eref{fepsnn} and \eref{momg} hold, and  if
$k_n\geq n$,  we have the upper bound
\begin{eqnarray}
{\mathbb E}\|g-\hat g_m^{(n)}\|^2\leq \|g-g_m\|^2+
\frac{m^2(M_2+1)}{n}+
\frac{2\lambda_1(f_\varepsilon,\kappa_0)\Gamma(m)}{n}+\frac{2R_{m}}{n}.\label{vit}
\end{eqnarray}
Finally, since $g_m$ is the orthogonal projection of $g$ on $S_m$,
we get that $g_m^*=g^*\mbox{1}\!\mbox{I}_{[-m\pi, m\pi]}$ and
therefore
$$
\|g-g_m\|^2 = \frac{1}{2 \pi} \|g^*-g_m^*\|^2 = \frac{1}{2 \pi}
\int_{|x|\geq \pi m} |g^*|^2(x)dx.
$$
If $g$ belongs to  the class $\mathcal{S}_{s,r,b}(C_1)$ defined in
\eref{super}, then
$$
\|g-g_m\|^2 \leq \frac{C_1}{2\pi} (m^2\pi^2+1)^{-s}\exp\{-2b
\pi^{r}m^{r}\}.
$$
Hence, according to \eref{vit}, if  \eref{momg} holds and $k_n\geq
n$, the risk of $\hat g_m^{(n)}$ is bounded by
\begin{multline*}
\frac{C_1}{2\pi}(m^2\pi^2+1)^{-s}\exp\{-2b\pi^{r}m^{r}\} +
\frac{2\lambda_1(f_\varepsilon,\kappa_0)(1+(\pi m)^2))^{\gamma}(\pi
m)^{1-\delta} \exp\left\{2\mu\pi^{\delta} m^{\delta}\right\}}{n} \\+
\frac{m^2(M_2+1)}{n} + \frac{2R_m}{n}.
\end{multline*}

Assume now that either  $\sum_{k >0} \beta_{{\bf X}, 1}(k) <
\infty$ or  $\sum_{k >0} \tau_{{\bf X}, 1}(k) < \infty$, so that
the residual terms $n^{-1}R_{m}+ n^{-1}m^2(M_2+1)$ are of order
$n^{-1}m^2$. As in the independent case, we choose $\breve m$ as
the minimizer of
$$
(m^2\pi^2+1)^{-s}\exp\{-2b\pi^{r}m^{r}\} + \frac{(\pi
m)^{2\gamma+1-\delta} \exp\left\{2\mu\pi^{\delta}
m^{\delta}\right\}}{n}.
$$
The behavior of $\breve m$ is recalled in Table~1. We see that in
all cases, the residual terms $n^{-1}R_{\breve m}+ n^{-1}{\breve
m}^2(M_2+1)$ of order $n^{-1}{\breve m}^2$ are negligible with
respect to the main terms since $n^{-1} \Delta(m)$ grows faster than $n^{-1}
m^2$ (recall that if $\delta=0$, we have the
restriction $\gamma >1/2$ (cf. Section \ref{add})). Hence the rate
of convergence of $\hat g^{(n)}_{\breve m}$  is the  same as in
the i.i.d. case (see Table~1 below).


%

\begin{table}[ptbh]
\caption{Choice of ${\breve m}$ and corresponding rates under
Assumptions \eref{condfeps}-\eref{fepsnn} and
\eref{super}.}\label{rates}
\begin{center}
{\small
\begin{tabular}{clcc}\cline{3-4}\cline{3-4}
\multicolumn{2}{c}{} &\multicolumn{2}{c}{$f_\varepsilon$} \\
\multicolumn{2}{c}{} & $\delta=0$ & $ \delta>0$ \\
\multicolumn{2}{c}{} & ordinary smooth & supersmooth \\\hline
\multirow{8}{.2cm}{\\\vfill\null $g$} & $\;$ & $\;$ & $\;$ \\
& $\begin{array}{l}
  r=0\\
  \small{\mbox{Sobolev}(s)}
\end{array}$ &
$\begin{array}{l}
  \pi {\breve m}=O(n^{1/(2s+2\gamma +1)})\\
  \mbox{rate}=O(n^{-2s/(2s+2\gamma+1)})\\
\mbox{{\it minimax rate}}
\end{array}$  &
$\begin{array}{l}
  \pi {\breve m}=[\ln(n)/(2\mu+1)]^{1/\delta}\\
  \mbox{rate}=O( (\ln(n))^{-2s/\delta})\\
\mbox{{\it minimax rate}}
\end{array}$ \\
\cline{2-4}
& $\begin{array}{l}
  r>0\\
  \mathcal{C}^\infty
\end{array}$ &
$\begin{array}{l} \\
  \pi {\breve m}=\left[{\ln(n)/2b}\right]^{1/r} \\
  \mbox{ rate}= \displaystyle  O\left(\frac{\ln(n)^{(2\gamma+1)/r}}n\right)\\
\mbox{{\it minimax rate}} \\
  \;\; \end{array}$ &
$\begin{array}{c}
  {\breve m}  \mbox{ solution of } \\
  {{\breve m}}^{2s+2\gamma+1-r}\exp\{2\mu
 (\pi {\breve{m}})^\delta+2b \pi^r {{\breve m}}^r\}\\
  \qquad= O(n)\\
\mbox{{\it minimax rate if }}r<\delta \mbox{ \it and } s=0
\end{array}$\\
\hline
\end{tabular}}
\end{center}
\end{table}

When $r>0, \delta>0$ the value of  ${\breve m}$ is not explicitly
given. It is obtained as the solution
 of the equation \begin{equation*}
{{\breve m}}^{2s+2\gamma+1-r}\exp\{2\mu
 (\pi {\breve{m}})^\delta+2b \pi^r {{\breve m}}^r\}= O(n).
 \end{equation*}
Consequently, the rate of $\hat g^{(n)}_{\breve m}$ is not explicit
and depends on the ratio $r/\delta$. If $r/\delta$ or $\delta/r$
belongs to  $ ]k/(k+1);(k+1)/(k+2)]$ with $k$ integer, the rate of
convergence can be expressed as a function of $k$. We refer to Comte
{\it et al.}~\citeyear{CRT1} for further discussions about those
rates. We refer to Lacour~\citeyear{CRASLAC} for explicit formulae
for the rates in  the special case $r>0,\delta>0$.

\section{Risk bounds for adaptive estimators} \label{sec5}
In the previous section, the construction of the estimators require
the knowledge of the smoothness of $g$. We now come to adaptive
estimation, without such  prior knowledge.

\subsection{A first bound in adaptive density deconvolution.}
 Theorem \ref{genepenmel} gives a general bound which holds under
mild dependence conditions, for $f_\varepsilon$ being either
ordinary or super smooth. For $a>1$, let pen$(m)$ be defined by
\begin{equation}\label{penalitess}
\displaystyle {\rm pen}(m)=\left\lbrace\begin{array}{l}
\displaystyle 24a\frac{\Delta(m)}{n}\,\mbox{
if } 0\leq \delta< 1/3,\\
\displaystyle
8a\left(1+\frac{48\mu\pi^\delta\lambda_2(f_\varepsilon,\kappa_0)}{\lambda_1(f_\varepsilon,\kappa_0')}\right)\frac{
\Delta(m)\, m^{\min((3\delta/2-1/2)_+,\delta))}}{n}\,\mbox{ if
}\delta \geq 1/3 .
\end{array}\right.
\end{equation}
The constant $\lambda_1(f_\varepsilon,\kappa_0)$ is defined in
\eref{lambda1} and $\lambda_2(f_\varepsilon,\kappa_0)$ is given by
\begin{eqnarray}\label{lambda2}
\lambda_2(f_\varepsilon,\kappa_0)=\parallel
f_\varepsilon\parallel\kappa_0^{-1}\sqrt{2\lambda_1(f_\varepsilon,\kappa_0)}\ind_{0\leq
\delta\leq 1}+2\lambda_1(f_\varepsilon,\kappa_0)\ind_{\delta>1}.
\end{eqnarray}
In order to bound up $\mathrm{pen}(m)$, we impose that
\begin{eqnarray}
\label{mn} \pi m_n\leq \left\{
\begin{array}{ll}
n^{1/(2\gamma+1)} &\mbox{ if }\delta=0\\
\displaystyle\left[\frac{\ln(n)}{2\mu}+\frac{2\gamma+1-\delta}{2\delta\mu
  }\ln\left(\frac{\ln(n)}{2\mu}\right)\right]^{1/ \delta} &\mbox{ if
} \delta>0.
\end{array}
\right.
\end{eqnarray}
Subsequently we set
\begin{eqnarray}
\label{kappaa}
\kappa_a=(a+1)/(a-1), \mbox{ and } C_a=\max(\kappa_a^2,2\kappa_a).\end{eqnarray}

\begin{thm}\label{genepenmel}
Assume that $f_\varepsilon$ satisfies
\eref{condfeps}-\eref{fepsnn},
 that $g$ satisfies \eref{momg}, and that $m_n$ satisfies \eref{mn}. Consider the collection of
estimators $\hat g_m^{(n)}$ defined by \eref{tronque} with $k_n\geq
n$ and $1\leq m\leq m_n$.
Let {\rm pen}$(m)$ be defined by \eref{penalitess}.
 The estimator $\tilde g = \hat g_{\hat m}^{(n)}$ defined by
\eref{estitronc} satisfies
\begin{equation*}
{\mathbb E}(\|g-\tilde g\|^2) \leq C_a\inf_{m\in \{1, \cdots,
m_n\}}\Big [\|g-g_m\|^2+ {\rm pen}(m)+\frac{m^2(M_2+1)}{n} \Big ] +
\frac{\overline{C}( R_{m_n}+ m_n)}{n}\, ,
\end{equation*} where $R_{m}$ is defined in (\ref{Rm}), $C_a$ is defined in \eref{kappaa}, and $\overline{C}$
is a constant depending on $f_{\varepsilon}$ and  $a$.
\end{thm}

 Let us compare the rate of $\tilde g$ with the rate obtained
in the independent framework. The term $\inf_{m\in
\{1,\cdots,m_n\}}[\|g-g_m\|^2+ {\rm pen}(m)+ m^2(M_2+1)/n]$
corresponds to the rate of $\tilde g$ when all variables are
i.i.d. The dependent context induces the additional term
$n^{-1}(R_{m_n}+m_n)$. If the dependence coefficients are summable
and the errors are super smooth, then $n^{-1}(R_{m_n}+m_n)$ is
negligible and $\tilde g$ achieves the rate of the independent
framework. If $\varepsilon$ is ordinary smooth, the term
$n^{-1}(R_{m_n}+m_n)$ may not be negligible and Theorem
\ref{genepenmel} is not precise enough.

\subsection{Adaptive density deconvolution for super smooth  $f_\varepsilon$.}

If  \eref{condfeps}-\eref{fepsnn} hold for some $\delta>0$, we have
the following corollary.
\begin{cor}\label{genepenmelss}
Assume that $f_\varepsilon$ satisfies
\eref{condfeps}-\eref{fepsnn} with $\delta>0$,  that $g$ satisfies
\eref{momg}, and that $m_n$ satisfies \eref{mn}.  Let $\mathrm
{pen}(m)$ be defined by (\ref{penalitess}). Consider the
collection of estimators $\hat g_m^{(n)}$ defined by
\eref{tronque} with $k_n\geq n$ and $1\leq m\leq m_n$.
\begin{enumerate}
\item If $\sum_{k>0} \beta_{{\bf X}, 1}(k) < \infty$, the
estimator $\tilde g = \hat g_{\hat m}^{(n)}$ defined by
\eref{estitronc} satisfies
\begin{eqnarray*}
{\mathbb E}(\|g-\tilde g\|^2) &\leq& C_a\inf_{m\in \{1, \cdots,
m_n\}}\Big[\|g-g_m\|^2+ {\rm pen}(m)+\frac{m^2(M_2+1)}{n} \Big ]
+ \frac{\overline{C}(\ln(n))^{1/\delta}}{n}
,
\end{eqnarray*}
where $C_a$ is defined in \eref{kappaa} and $\overline{C}$ is a
constant depending on $f_\varepsilon$, $a$ and $\sum_{k>0} \beta_{{\bf X}, 1}(k)$.
\item If   $\sum_{k>0} \tau_{{\bf X}, 1}(k) < \infty$,  the
estimator $\tilde g = \hat g_{\hat m}^{(n)}$ defined by
\eref{estitronc} satisfies
\begin{eqnarray*}
{\mathbb E}(\|g-\tilde g\|^2) &\leq& C_a\inf_{m\in \{1, \cdots,
m_n\}}\Big[\|g-g_m\|^2+ {\rm pen}(m)+\frac{m^2(M_2+1)}{n} \Big ]
+\frac{\overline{C}(\ln(n))^{2/\delta}}{n}
,
\end{eqnarray*}
where $C_a$ is defined  in \eref{kappaa} and $\overline{C}$ is a
constant depending on $f_\varepsilon$, $a$ and  $\sum_{k>0} \tau_{{\bf X}, 1}(k)$.
\end{enumerate}
\end{cor}

Corollary \ref{genepenmelss} requires important comments. The
terms involving  power of $\ln(n)$ are negligible with respect
to $\inf_{m\in \{1, \cdots, m_n\}}[\|g-g_m\|^2+ {\rm
pen}(m)+m^2(M_2+1)/n]$. The risk of $\tilde g$ is of order
$\inf_{m\in \{1,\cdots,m_n\}}[\|g-g_m\|^2+\mbox{pen}(m)]$, that is
of the best order,  as in the independent framework. The penalty
does not depend on the dependence coefficients and is the same as in
the independent framework.

As a conclusion, we see that the adaptive estimator $\tilde g$ built
with the same penalty as in the independent framework, still
achieves the best rates under  mild conditions on the dependence
coefficients.

\subsection{Adaptive density deconvolution for ordinary smooth  $f_\varepsilon$.}

For $a>1$, define {\rm pen}$(m)$ by
\begin{equation}\label{penaliteos}
{\rm pen}(m)= \frac{25a\Delta(m)}{n} \, .
\end{equation}

\begin{thm}\label{genepenmelos}
Assume that $f_\varepsilon$ satisfies
\eref{condfeps}-\eref{fepsnn} with $\delta=0$, that $g$ satisfies
\eref{momg}, and that $m_n$ satisfies \eref{mn}. Let
$\mathrm{pen}(m)$ be defined by \eref{penaliteos}. Consider the
collection of estimators $\hat g_m^{(n)}$ defined by
\eref{tronque} with $k_n\geq n$ and $1\leq m\leq m_n$.
\begin{enumerate}
\item If  $\beta_{\mathbf{X},\infty}(k)=O(k^{-(1+ \theta)})$ for
some $\theta
>(2\gamma+3)/(2\gamma+1)$, then the  estimator  $\tilde g = \hat
g_{\hat m}^{(n)}$ defined by \eref{estitronc} satisfies
\begin{equation}\label{resuosbeta}
{\mathbb E}(\|g-\tilde g\|^2) \leq C_a\inf_{m\in \{1, \cdots, m_n
\}}\Big [\|g-g_m\|^2+ {\rm pen}(m)+\frac{m^2(M_2+1)}{n} \Big ] +
\frac{\overline{C}}{n},
\end{equation} where $C_a$ is defined in \eref{kappaa} and
 $\overline{C}$ is a constant depending on $f_{\varepsilon}$, $a$,
and  $\sum_{k>0}\beta_{\mathbf{X},\infty}(k)$.
\item If $\tau_{\mathbf{X},\infty}(k)= O(k^{-(1+ \theta)})$ for
some $\theta
>(2\gamma+5)/(2\gamma+1)$, then  the estimator $\tilde g = \hat
g_{\hat m}^{(n)}$ defined by \eref{estitronc} satisfies
(\ref{resuosbeta}),
where $\overline{C}$  is a constant depending on $f_{\varepsilon}$,
$a$ and $\sum_{k
>0}\tau_{\mathbf{X},\infty}(k)$.
\end{enumerate}
\end{thm}

\begin{rem}{\rm
Note that the condition for $\beta_{\mathbf{X},\infty}(k)$ is
realized for any $\gamma>1/2$ provided $\theta >2$.  In the same
way, the condition for $\tau_{\mathbf{X},\infty}(k)$ is realized for
any $\gamma >1/2$ provided $\theta >3$. In both cases, the condition
on $\theta$  is weaker as $\gamma$ increases. In other words, the
smoother is $f_\varepsilon$, the weaker is the condition on the
dependence coefficients. }\end{rem}

\begin{rem}
\label{remos}{\rm For $m$ large enough, the penalty function given
for ordinary smooth errors in Theorem \ref{genepenmelos} is an upper
bound of more precise penalty functions which depend on the
dependence coefficients. Under the assumptions of (1) in Theorem
\ref{genepenmelos},  let {\rm pen}$(m)$ be defined by
 \begin{equation}\label{pencorbeta}
{\rm pen}(m)= \frac{24a\Delta(m) + 128a\Big(1 +
4\sum_{k=1}^{n} \beta_{\mathbf{X},1}(k)\Big){m}}{n}.
\end{equation}

Under the assumptions of (2) in Theorem \ref{genepenmelos} let {\rm
pen}$(m)$ be defined by
 \begin{equation}\label{pencortau}
{\rm pen}(m)=\frac{24a\Delta(m)}{n}+\frac{64a\left[1+38\ln(m)\right]
\left(m+\pi\sum_{k=1}^{n}\tau_{\mathbf{X},1}(k)m^2\right)}{n}
\end{equation}
In both cases,  the estimator $\tilde g = \hat g_{\hat m}^{(n)}$
defined by \eref{estitronc} satisfies (\ref{resuosbeta}).
Remark \ref{remos} follows from the proof of Theorem \ref{genepenmelos}.

}\end{rem}

\subsection{Case without noise}

One can deduce from Proposition \ref{Vitss}, Theorem \ref{genepenmelos}, its
proof and Remark \ref{remos}, a result for density estimation without errors, on the whole real
line, that is when the $X_i$ is observed. If $\varepsilon=0$, then we can consider that $Z=X$ and replace $f_{\varepsilon}^*$ by 1.
It follows that $u_t^*(Z_i)=t(X_i)$ and the contrast $\gamma_n$ simply becomes
\begin{equation}\label{density}
\gamma_{n,X}(t)=\|t\|^2-\frac 2n\sum_{i=1}^n t(X_i).
\end{equation}
Let $k_n \geq n^2$, and consider as previously
\begin{equation}\label{gamnX} \hat g_m^{(n)} =\arg\min_{t\in S_m^{(n)}}
\gamma_{n,X}(t), \;\; {\rm pen}(m)= 128a\Big(1+4
\sum_{k=1}^n\beta_{{\bf X}, 1}(k)\Big)\frac mn \, ,
\end{equation}
and 
\begin{equation}
\label{gamnX1}
\hat m=\arg\min_{m\in \{1, \dots, n\}} [\gamma_{n,X}(g_m^{(n)}) + {\rm
pen}(m)].
\end{equation}
The following results follow straightforwardly.

\begin{cor}\label{densitycor}
Assume that $\varepsilon=0$. Let $k_n \geq n^2$. Then
\begin{enumerate}
\item 
$$
 {\mathbb E}\|g-\hat g_m^{(n)}\|^2\leq \|g-g_m\|^2+
\frac{m(M_2+3)}{n} +  \frac{2R_{m}}{n}.
$$

\item  If $\beta_{{\bf X},\infty}=O(k^{-(1+\theta)})$ for
some $\theta>3$, then the estimator $\tilde g=\hat g_{\hat
m}$ defined by \eref{gamnX} and \eref{gamnX1} satisfies
\begin{equation*}
{\mathbb E}(\|g-\tilde g\|^2) \leq C_a\inf_{m\in \{1, \cdots, n
\}}\Big [\|g-g_m\|^2+ {\rm pen}(m)+\frac{m(M_2+1)}{n} \Big ] +
\frac{\overline{C}}{n},
\end{equation*} where $C_a$ is defined in \eref{kappaa} and
 $\overline{C}$ is a constant depending on  $a$
and  $\sum_{k>0}\beta_{\mathbf{X},\infty}(k)$.
\end{enumerate}
\end{cor}

The result  (1) shows that if $\sum_{k >0} \beta_{{\bf X}, 1}(k) < \infty$, one obtains
the same bounds (and the same rates) as in the i.i.d. case. However,
if $\sum_{k >0} \tau_{{\bf X}, 1}(k) < \infty$ the term $n^{-1} R_m$
is of order $n^{-1} m^2$ and the rates for $\hat g_m^{(n)}$ are less good than in the
i.i.d. case. 
 
The result (2)  shows that this estimation procedure also works in
density estimation without errors. It allows to estimate a density
on the whole real line and to reach the usual rates of convergence,
by using a penalty of the classical order $m/n$. This remark is
valid in the $\beta$-mixing framework and  in the case
 of independent $X_i$'s. We refer to Pensky~\citeyear{Pensky99} and Rigollet~\citeyear{rigollet} for recent results
 in
adaptive density estimation on the whole real line in the i.i.d.
case.

\section{Proofs}

\setcounter{equation}{0}
\setcounter{lem}{0}
\setcounter{thm}{0}
\subsection{Proof of Proposition \ref{Vitss}}
The proof of the proposition \ref{Vitss} follows the same lines as
in the independent framework (see Comte  {\it et
al.}~\citeyear{CRT1}). The main difference lies in the control of
the variance term. We keep the same notations as in Section
\ref{nonpen}. According to \eref{tronque}, for any given $ m$
belonging to $\{1,\cdots,m_n\}$, $\hat g_{m}^{(n)}$ satisfies,
$\gamma_n(\hat g_{m}^{(n)} )-\gamma_n(g_{m}^{(n)})\leq 0.$ For a
random variable $Y$ with density $f_Y$, and any function $\psi$ such
that $\psi(Y)$ is integrable, let
\begin{eqnarray}\label{nu}
\nu_{n,Y}(\psi)=\frac 1n \sum_{i=1}^n [\psi(Y_i)-\langle \psi,
f_Y\rangle] ,\mbox{ so that }\;\;\;\nu_{n,Z}(u_t^*)=\frac 1n
\sum_{i=1}^n \left[u_t^*(Z_i)-\langle t,g\rangle
\right].\end{eqnarray} Since
\begin{eqnarray}
\label{difgamma}
\gamma_n(t)-\gamma_n(s)=\|t-g\|^2-\|s-g\|^2-2\nu_{n,Z}(u_{t-s}^*),
\end{eqnarray}
we infer that \begin{equation}\label{premineq}\|g-\hat
g_m^{(n)}\|^2\leq \|g-g_m^{(n)}\|^2 + 2\nu_{n,Z}\left(u_{\hat
g_m^{(n)} - g_m^{(n)}}^*\right) \, . \end{equation}
 Writing that $\hat
a_{m,j}-a_{m,j}= \nu_{n,Z}(u_{\varphi_{m,j}}^*)$, we obtain
$$
\nu_{n,Z}\left(u_{\hat g_m^{(n)}-g_m^{(n)}}^*\right)=\sum_{\vert
j\vert \leq k_n} (\hat
a_{m,j}-a_{m,j})\nu_{n,Z}(u_{\varphi_{m,j}}^*) = \sum_{\vert
j\vert\leq k_n} [\nu_{n,Z}(u_{\varphi_{m,j}}^*)]^2.$$
Consequently,
$\esp\|g-\hat g_m^{(n)}\|^2\leq  \|g-g_m^{(n)}\|^2 +
2\sum_{j\in \mathbb{Z}}\mathbb{E}[(\nu_{n,Z}(u_{\varphi_{m,j}}^*))^2].$
According to Comte {\it et al.} \citeyear{CRT1},
\begin{equation} \label{biasineq}
\|g-g_m^{(n)}\|^2= \parallel
g-g_m\parallel^2+\|g_m-g_m^{(n)}\|^2\leq \parallel g-g_m\parallel^2+
\frac{(\pi m)^2(M_2+1)}{k_n}.
\end{equation}
The variance term is studied by using that for $f\in
\mathbb{L}_1(\mathbb{R})$,
\begin{eqnarray}
\label{lemnu} \nu_{n,Z}(f^*)=\int \nu_{n,Z}(e^{ix\cdot})f(x)dx.
\end{eqnarray}
Now, we use \eref{lemnu}
and apply Parseval's formula to obtain
\begin{equation}\label{varineq}
\mathbb{E}\left(\sum_{j\in
\mathbb{Z}}(\nu_{n,Z}(u_{\varphi_{m,j}}^*))^2\right)=\frac{1}{4\pi^2}\sum_{j\in
\mathbb{Z}}\mathbb{E}\Big(\int
\frac{\varphi_{m,j}^*(-x)}{f_\varepsilon^*(x)}\nu_{n,Z}(e^{ix\cdot})dx\Big)^2=\frac{1}{2\pi}\int_{-\pi
m}^{\pi m}\frac{\mathbb{E}\vert
\nu_{n,Z}(e^{ix\cdot})\vert^2}{\vert f^*_\varepsilon(x)\vert^2}dx.
\end{equation}
Since $\nu_{n,Z}$ involves centered and stationary variables,
\begin{eqnarray}
\mathbb{E}\vert
\nu_{n,Z}(e^{ix\cdot})\vert^2&=&\mbox{Var}\vert
\nu_{n,Z}(e^{ix\cdot})\vert=\frac{1}{n^2}\left(\sum_{k=1}^n\mbox{Var}(e^{ixZ_k})+\sum_{1\leq
k\not=l\leq
n}\mbox{Cov}(e^{ixZ_k},e^{ixZ_l})\right)\nonumber \\
&=&\frac{1}{n}\mbox{Var}(e^{ixZ_1})+\frac{1}{n^2}\sum_{1\leq
k\not=l\leq n}\mbox{Cov}(e^{ixZ_k},e^{ixZ_l}). \label{devpt}
\end{eqnarray}
Since $(X_i)_{i \geq 1}$ and $(\varepsilon_i)_{i \geq 1}$ are
independent, we have $\mathbb{E}(e^{ix
Z_k})=f_\varepsilon^*(x)g^*(x)$ so that
\begin{eqnarray*}
\mbox{Cov}(e^{ixZ_k},e^{ixZ_l})=\mathbb{E}(e^{ix(Z_l-Z_k)})-\vert\mathbb{E}(e^{ixZ_k})\vert^2=\mathbb{E}(e^{ix(Z_l-Z_k)})-\vert
f_\varepsilon^*(x)g^*(x)\vert^2.
\end{eqnarray*}
Next, by independence of $X$ and $\varepsilon$, we write, for
$k\neq l$,
\begin{eqnarray*}
\mathbb{E}(e^{ix(Z_l-Z_k)})=\mathbb{E}(e^{ix(X_l-X_k)}){\mathbb
E}(e^{ix(\varepsilon_l-\varepsilon_k)})={\mathbb
E}(e^{ix(X_l-X_k)})|f_\varepsilon^*(x)\vert^2,
\end{eqnarray*}
and consequently
\begin{eqnarray}\label{liencov}
\mbox{Cov}(e^{ixZ_k},e^{ixZ_l})=\mbox{Cov}(e^{ixX_k},e^{ixX_l})\vert
f_{\varepsilon}^*(x)\vert^2.\end{eqnarray} From \eref{devpt},
  \eref{liencov} and the stationarity of $(X_i)_{i \geq 1}$, we obtain
  that
\begin{equation}\label{nunfin}
\mathbb{E}\vert \nu_{n,Z}(e^{ix\cdot})\vert^2 \leq
\frac{1}{n}+\frac{2}{n}\sum_{k=2}^n\left\vert\mbox{Cov}(e^{ixX_1},e^{ixX_k})\right\vert\vert
f_{\varepsilon}^*(x)\vert^2.
\end{equation}
The first part of Proposition \ref{Vitss} follows from  the
stationarity of the $X_i$'s, and from \eref{premineq},
\eref{biasineq}, \eref{varineq} and \eref{nunfin}.

Let us prove that $R_m \leq \min ( R_{m, \beta}, R_{m, \tau})$,
where  $R_{m, \beta}$ and $R_{m, \tau}$ are defined in Proposition
\ref{Vitss}. Using  the inequalities \eref{ibrb} and \eref{evident},
we obtain the  bounds
\begin{eqnarray*}
    |\mbox{Cov}(e^{ixX_1},e^{ixX_k})| \leq 2 \beta_{{\bf X}, 1} (k-1)
    \quad \text{ and }
    |\mbox{Cov}(e^{ixX_1},e^{ixX_k})| &\leq&   |x| \tau_{{\bf X}, 1}
   (k-1)
\end{eqnarray*}
(for the last inequality, note that $t \rightarrow e^{ixt}$ is
$|x|$-Lipschitz).  The result easily follows.

\subsection{Proof of Theorem \ref{genepenmel}}
By definition, $\tilde g$ satisfies that for all $m\in
\{1,\cdots,m_n\}$, $$\gamma_n(\tilde g)+\mbox{pen}(\hat m)\leq
\gamma_n(g_m)+\mbox{pen}(m).$$ Therefore, by using \eref{difgamma}
we get that
\begin{eqnarray*}
\|\tilde g-g\|^2\leq \|g_m^{(n)}-g\|^2+2\nu_{n,Z}(u_{\tilde
g-g_m^{(n)}}^*)+\mbox{pen}(m)-\mbox{pen}(\hat m).
\end{eqnarray*}
If $t=t_1+t_2$ with $t_1$ in $S_m^{(n)}$ and $t_2$ in
$S_{m'}^{(n)}$, $t^*$ has its support in $[-\pi {\max(m,m')}, \pi
{\max(m,m')}]$ and  $t$ belongs to $S_{\max(m,m')}^{(n)}$. Set
 $B_{m, m'}(0,1)=\{t\in
S_{\max(m,m')}^{(n)} \;/\; \|t\|=1\}.$ For $\nu_{n,Z}$ defined in \eref{nu} we
get $$|\nu_{n,Z}(u_{\tilde
g-g_m^{(n)}}^*) |\leq \|\tilde g-g_m^{(n)}\|\sup_{t\in
  B_{m,\hat m}(0,1)}|\nu_{n,Z}(u_t^*)|.$$
Using that $2uv \leq a^{-1}u^2+av^2$ for any $a>1$, leads to
\begin{eqnarray*} \|\tilde g-g\|^2 &\leq& \|g_m^{(n)}
-g\|^2 + a^{-1}\|\tilde g-g_m^{(n)}\|^2  + a\sup_{t\in B_{m,\hat
m}(0,1)}(\nu_{n,Z}(u_t^*))^2+ {\rm pen}(m)- {\rm pen}(\hat m).
\end{eqnarray*}
Now, according to Lemma \ref{nucond}, write that $\nu_{n,Z}(u_t^*)
=\nu_n^{(1)}(t)+\nu_{n,X}(t),$ where
\begin{eqnarray}
\label{nu12}
\nu_n^{(1)}(t)=n^{-1}\sum_{i=1}^n[ u_t^*(Z_i)-\mathbb{E}(u_t^*(Z_i)|\sigma(X_i,\,i\geq
1))]=n^{-1}\sum_{i=1}^n [u_t^*(Z_i)-t(X_i)].
\end{eqnarray}
Consequently,
\begin{eqnarray*} \|\tilde g-g\|^2 &\leq &\|g_m^{(n)}
-g\|^2 + a^{-1}\|\tilde g-g_m^{(n)}\|^2  + 2a\sup_{t\in B_{m,\hat
m}(0,1)}(\nu_n^{(1)}(t))^2+ 2a\sup_{t\in B_{m,\hat
m}(0,1)}(\nu_{n,X}(t))^2\\&&+{\rm pen}(m)- {\rm pen}(\hat m).
\end{eqnarray*}
Hence by writing that $\|\tilde g-g_m^{(n)}\|^2\leq (1+\kappa_a^{-1})\|\tilde
g-g\|^2+(1+\kappa_a)\|g-g_m^{(n)}\|^2$ with $\kappa_a$ defined in \eref{kappaa}, we have
\begin{eqnarray*}
\|\tilde g-g\|^2 &\leq& \kappa_a^2\|g_m^{(n)} -g\|^2 +
2a\kappa_a\sup_{t\in B_{m,\hat m}(0,1)}(\nu_n^{(1)}(t))^2+
2a\kappa_a\sup_{t\in B_{m,\hat
m}(0,1)}(\nu_{n,X}(t))^2\\&&+\kappa_a({\rm pen}(m)- {\rm pen}(\hat
m)).
\end{eqnarray*}
Choose some positive function $p(m,m')$ such that
\begin{eqnarray}
\label{pmmp} 2ap(m,m')\leq \mbox{pen}(m)+\mbox{pen}(m').
\end{eqnarray} For this function $p(m,m')$  we have
\begin{eqnarray}
\|\tilde g-g\|^2&\leq&
\kappa_a^2\|g-g_m^{(n)}\|^2+2\kappa_a\mbox{pen}(m)+2a\kappa_a\sup_{t\in
B_{m,\hat m}(0,1)}(\nu_{n,X}(t))^2+2a\kappa_a W_n(m,\hat
m)\nonumber\\&\leq&\label{majo2}
\kappa_a^2\|g-g_m^{(n)}\|^2+2\kappa_a\mbox{pen}(m)+2a\kappa_a\sup_{t\in
B_{m,\hat m}(0,1)}(\nu_{n,X}(t))^2+2a\kappa_a
\sum_{m'=1}^{m_n}W_n(m,m'),
\end{eqnarray}
where
\begin{equation}\label{Wg} W_n(m,m'):=\Big[\sup_{t\in B_{m, m'}(0,1)}
 |\nu_n^{(1)}(t)|^2-p(m, m')\Big]_+,
\end{equation}
The main parts of the proof lies in the two following points~:

\textbf{1)} Study of $W_n(m,m')$, and more precisely find
$p(m,m')$ such that for a constant $A_1$,
\begin{equation}\label{but1}
\sum_{m'=1}^{m_n} \mathbb{E}(W_n(m,m'))\leq \frac{A_1}{n}.
\end{equation}

\textbf{2)} Study of $\sup_{t\in B_{m,\hat
m}(0,1)}(\nu_{n,X}(t))^2$ and more precisely prove that
\begin{eqnarray}
\label{but2} \mathbb{E}\Big[\sup_{t\in B_{m,\hat
m}(0,1)}(\nu_{n,X}(t))^2\Big]\leq \frac{{m_n}+R_{m_n}}{n},
\end{eqnarray}
where $R_{m}$ is defined in \eref{Rm}. Combining \eref{majo2},
\eref{but1} and \eref{but2}, we infer that, for all $1 \leq m \leq
m_n$
$$  \mathbb{E}\|g-\tilde g\|^2 \leq \kappa_a^2
  \|g-g_m^{(n)}\|^2
  +2\kappa_a {\rm pen}(m)+\frac{2a\kappa_a({m_n}+R_{m_n})}{n}
+\frac{2a\kappa_a A_1}{n}.$$ If we denote by $C_a=\max(\kappa_a^2,
2\kappa_a)$, this can also be written
\begin{eqnarray*}
\mathbb{E}\|g-\tilde g\|^2  \!\!\!&\leq& \!\!\!C_a\inf_{m\in
\{1,\cdots,m_n\}} \big[ \|g-g^{(n)}_m\|^2 +\|g_m^{(n)}-g_m\|+ {\rm
pen}(m)\big]+ \frac{2a\kappa_a(L_{m_n}+R_{m_n})}{n}  + \frac{2a\kappa_a A_1}{n}\\
\!\!\!&\leq & \!\!\!
 C_a\inf_{m\in
\{1,\cdots,m_n\}} \big[ \|g-g_m\|^2 +(M_2+1)m^2/k_n+ {\rm
pen}(m)\big]+ \frac{2a\kappa_a(L_{m_n}+R_{m_n})}{n} +
\frac{2a\kappa_a A_1}{n}.
\end{eqnarray*}

\noindent \textbf{Proof of \eref{but1}} We start by writing
$\mathbb{E}(W_n(m,m'))=\mathbb{E}[\sup_{t\in B_{m, m'}(0,1)}
 |\nu_n^{(1)}(t)|^2-p(m, m')]_+$ as
$$\mathbb{E}\Big \{\mathbb{E}_{\textbf{X}}\Big[\sup_{t\in B_{m, m'}(0,1)}
 |\nu_n^{(1)}(t)|^2-p(m, m')\Big]_+\Big\},$$
where $\mathbb{E}_{\textbf{X}}(Y)$ denotes the conditional
expectation $\mathbb{E}(Y|\sigma(X_i,\,i\geq 0)).$ The point is
that, conditionally to $\sigma(X_i,\,i\geq 0)$, the random variables
$u_t^*(Z_i)-\mathbb{E}(u_t^*(Z_i)|\sigma(X_i,\,i\geq 0))$ are
centered, independent but non identically distributed. We proceed as
in the independent case (see Comte {\it et al.}~\citeyear{CRT1}), by
applying the following  Lemma  to the expectation
$\mathbb{E}_{\textbf{X}}[\sup_{t\in B_{m, m'}(0,1)}
 |\nu_n^{(1)}(t)|^2-p(m, m')]_+$.
\begin{lem}
\label{Concent} Let $Y_1, \dots, Y_n$ be independent random
variables and  let ${\mathcal F}$ be a countable class of uniformly
bounded measurable functions. Then for $\xi^2>0$
\begin{equation*}
 \mathbb{E}\Big[\sup_{f\in {\mathcal
F}}|\nu_{n,Y}(f)|^2-2(1+2\xi^2)H^2\Big]_+ \leq \frac
4{K_1}\left(\frac vn e^{-K_1\xi^2 \frac{nH^2}v} +
\frac{98M_1^2}{K_1n^2C^2(\xi^2)} e^{-\frac{2K_1
C(\xi)\xi}{7\sqrt{2}}\frac{nH}{M_1}}\right),
\end{equation*}
with $C(\xi)=\sqrt{1+\xi^2}-1$, $K_1=1/6$, and
$$\sup_{f\in {\mathcal F}}\|f\|_{\infty}\leq M_1, \;\;\;\;
\mathbb{E}\Big[\sup_{f\in {\mathcal F}}|\nu_{n,Y}(f)|\Big]\leq H,
\;\;\;\; \sup_{f\in {\mathcal F}}\frac{1}{n}\sum_{k=1}^n{\rm
Var}(f(Y_k)) \leq v.$$\end{lem} The proof of this inequality can be
found in Appendix. It comes from  a concentration Inequality in
Klein and Rio \citeyear{KleinRio}
 and arguments
that can be found in Birg\'e and Massart \citeyear{BirgeMassart98}.
Usual density arguments show that this result can be applied to the class of
functions ${\mathcal F}= B_{m,m'}(0,1)$.
 Let us denote by $m^*=\max(m,m')$. Applying Lemma
 \ref{Concent}, one has the bound
\begin{equation*}
 \mathbb{E}_{\textbf{X}}\Big[\sup_{t\in B_{m,m'}(0,1)}|\nu_{n}^{(1)}(t)|^2-2(1+2\xi^2)H^2\Big]_+ \leq \frac
6{K_1}\left(\frac vn e^{-K_1\xi^2 \frac{nH^2}v} +
\frac{98M_1^2}{K_1n^2C^2(\xi^2)} e^{-\frac{K_1
C(\xi)\xi}{7\sqrt{2}}\frac{nH}{M_1}}\right),
\end{equation*}
where
$$\sup_{t \in B_{m,m'}(0,1)}\|u_t^*(Z_1)\|_{\infty}\leq M_1, \;\;\;\;
\mathbb{E}_{\textbf{X}}\Big[\sup_{t \in
B_{m,m'}(0,1)}|\nu_{n}^{(1)}(t)|\Big]\leq H, \;\;\;\; \sup_{t\in
B_{m,m'}}\frac{1}{n}\sum_{k=1}^n{\rm Var}_{\textbf{X}}(u_t^*(Z_k))
\leq v.$$ By applying Lemma \ref{l3}, we propose to take
$$H^2= H^2(m^*)= \frac{\Delta(m^*)}{n},
\quad M_1=M_1(m^*) =\sqrt{nH^2}\mbox{ and
}v=v(m^*)=\frac{\sqrt{\Delta_2(m^*,h)}}{2\pi}$$ with, for $f_Z$
denoting the density of $Z_1$,
\begin{equation}\label{Delta2}
\Delta_2(m, h)=  \int_{-\pi m}^{\pi m} \int_{-\pi m}^{\pi
m}\frac{\vert f_Z^*(x-y)\vert^2}{\vert
f_\varepsilon^*(x)f_\varepsilon^*(y)\vert^2}dxdy.
 \end{equation}
 From the definition \eref{Wg} of $W_n(m,m')$, by taking
$p(m,m')=2(1+2\xi^2)H^2(m^*)$, we get that
\begin{eqnarray}
\mathbb{E}(W_n(m,m'))\leq\mathbb{E}\Big\{\mathbb{E}_{\textbf{X}}\Big[\sup_{t\in
B_{m,m'}(0,1)}
  |\nu_n^{(1)}(t)|^2
  -2(1+2\xi^2)H^2(m^*)\Big]_+\Big\}.\label{bWn}\end{eqnarray}
  According to the condition \eref{pmmp}, we thus take
  pen$(m)=4ap(m,m)=8n^{-1}a(1+2\xi^2)\Delta(m)$ where $\xi^2$ is
  suitably chosen in the control of the sum of the right-hand side
  of \eref{bWn}.
Set $m_0$ such that for $m^*\geq m_0$
\begin{eqnarray}
\label{encadDelta}
(1/2)\lambda_1(f_\varepsilon,\kappa_0')\Gamma(m^*)\leq\Delta(m^*)\leq
2\lambda_1(f_\varepsilon,\kappa_0)\Gamma(m^*)
\end{eqnarray}
where $\Gamma(m)$ is defined in \eref{gammadem} and
$\lambda_1(f_\varepsilon,\kappa_0)$  and
$\lambda_1(f_\varepsilon,\kappa_0')$ are defined in \eref{lambda1}.
We split the sum over $m'$ in two parts and write
\begin{eqnarray}
\sum_{m'=1}^{m_n}\mathbb{E}(W_n(m,m'))=\sum_{m'|m^*<
m_0}\mathbb{E}(W_n(m,m'))+ \sum_{m'|m^*\geq
m_0}\mathbb{E}(W_n(m,m'))\label{m0}.\end{eqnarray} By applying
Lemma \ref{Concent} and \eref{encadDelta}, we get the global bound
$\mathbb{E}_{\textbf{X}}(W_n(m,m'))\leq K[I({m^*})+II(m^*)],$
where $I(m^*)$ and $II(m^*)$ are defined by
\begin{eqnarray*}
&&I(m^*)=\frac{v(m^*)}{n}
\exp\left\lbrace-K_1\xi^2\frac{\Delta(m^*)}{v(m^*)}\right\rbrace\\
\mbox{ and } &&II(m^*)= \frac{\Delta(m^*)}{n^2}
\exp\left\lbrace-\frac{2K_1\xi
C(\xi)}{7\sqrt{2}}\sqrt{n}\right\rbrace,
\end{eqnarray*}
Since  $I$ and $II$ do not depend on the $X_i$'s,  we infer that
$\mathbb{E}(W_n(m,m'))\leq K[I({m^*})+II(m^*)].$

When $m^*\leq m_0$,
with $m_0$ finite, we get that for all $m\in \{1,\cdots,m_n\}$,
\begin{eqnarray*}
\sum_{m'|m^*\leq m_0}\mathbb{E}(W_n(m,m'))\leq \frac{C(m_0)}{n}.
\end{eqnarray*}
We now come to the sum over $m'$ such that $m^*> m_0$.

When $\delta>1$ we use a rough bound for $\Delta_2(m,h)$ given by
 $\sqrt{\Delta_2(m,h)}\leq 2\pi nH^2(m)$.

When $0\leq\delta \leq 1$, write that
\begin{eqnarray*}
\Delta_2(m,h)\leq \parallel
|f_\varepsilon^*|^{-2}\ind_{[-\pi m,\pi m]}\parallel_\infty
\Delta(m)\parallel h^*\parallel^2(2\pi).
\end{eqnarray*}
Under \eref{condfeps}-\eref{fepsnn} we use that $\|h^*\|^2 \leq
\|f_\varepsilon^*\|^2<\infty$, that
$\sqrt{2\pi}\|f_\varepsilon^*\|=\|f_\varepsilon\|$ and apply
\eref{encadDelta} to infer that for $m^*\geq m_0$,
\begin{eqnarray}\label{vstar} v(m^*)=\frac{\sqrt{\Delta_2(m^*,h)}}{2\pi}\leq
\lambda_2(f_\varepsilon,\kappa_0) \Gamma_2(m^*),
\end{eqnarray}
where $\lambda_2(f_\varepsilon,\kappa_0)$ is defined in
\eref{lambda2}
 and
\begin{eqnarray} \label{Gamma2} \quad\Gamma_2(m)=(1+(\pi
m)^2)^{\gamma}(\pi m)^{\min((1/2-\delta/2),(1-\delta))}\exp(2\mu(\pi
m)^\delta)=(\pi m)^{-(1/2-\delta/2)_+}\Gamma(m).
\end{eqnarray}
Combining \eref{encadDelta} and \eref{vstar}, we get that for $m^*\geq m_0$,
\begin{eqnarray*}
&&I(m^*)\leq\frac{\lambda_2(f_\varepsilon,\kappa_0)\Gamma_2(m^*)}{n}
\exp\left\lbrace-\frac{K_1\xi^2\lambda_1(f_\varepsilon,\kappa_0')}{2\lambda_2(f_\varepsilon,\kappa_0)}(\pi m^*)^{(1/2-\delta/2)_+}\right\rbrace\\
\mbox{ and } &&II(m^*)\leq\frac{\Delta(m^*)}{n^2}
\exp\left\lbrace-\frac{2K_1\xi
C(\xi)\sqrt{n}}{7\sqrt{2}}\right\rbrace.
\end{eqnarray*}

\noindent $\bullet$ Study of $\sum_{m'|m^*\geq m_0}II(m^*)$.
According to the choices for $v(m^*)$, $H^2(m^*)$ and $M_1(m^*)$, we
have
\begin{eqnarray*}
\sum_{m'|m^*\geq m_0}II(m^*)&\leq&
\sum_{m'=1}^{m_n}\frac{\Delta(m^*)}{n^2}
\exp\left\lbrace\frac{-2K_1\xi
C(\xi)\sqrt{n}}{7\sqrt{2}}\right\rbrace\\&\leq&
\frac{\Delta(m_n)m_n}{n^2} \exp\left\lbrace\frac{-2K_1\xi C(\xi)
  \sqrt{n}}{7\sqrt{2}}\right\rbrace.
\end{eqnarray*}
Since under \eref{mn}, $n^{-1}\Delta(m_n)$ is bounded, we deduce
that $\sum_{m'|m^*\geq m_0}II(m^*)\leq n^{-1}C$.

\medskip

\noindent $\bullet$ Study of $\sum_{m'|m^*\geq m_0}I(m^*)$. Denote
by $\psi=2\gamma+ \min(1/2-\delta/2,1-\delta)$,
$\omega=(1/2-\delta/2)_+$, and
$K'=K_1\lambda_1(f_\varepsilon,\kappa_0')/(2\lambda_2(f_\varepsilon,\kappa_0))$.
For $a,b\geq 1$, we have that
\begin{eqnarray}\nonumber
\max(a,b)^{\psi}e^{2\mu\pi^{\delta}
\max(a,b)^{\delta}}e^{-K'\xi^2\max(a,b)^{\omega}} &\leq&
(a^{\psi}e^{2\mu\pi^{\delta} a^{\delta}}+b^{\psi}
e^{2\mu\pi^{\delta}
b^{\delta}})e^{-(K'\xi^2/2)(a^{\omega} + b^{\omega})}\\
&\leq&   a^{\psi}e^{2\mu\pi^{\delta}
a^{\delta}}e^{-(K'\xi^2/2)a^{\omega}} e^{-(K'\xi^2/2) b^{\omega}} +
b^{\psi} e^{2\mu\pi^{\delta} b^{\delta}} e^{-(K'\xi^2/2)
b^{\omega}}.\label{eqmax}\end{eqnarray} Consequently,
\begin{eqnarray}
\sum_{m'|m^*\geq m_0}I(m^*)&\leq&\!
\sum_{m'=1}^{m_n}\frac{\lambda_2(f_\varepsilon,\kappa_0)\Gamma_2(m^*)}{n}
\exp\left\lbrace-\frac{K_1\xi^2(\lambda_1(f_\varepsilon,\kappa_0')}{2\lambda_2(f_\varepsilon,\kappa_0)}(\pi
m^*)^{(1/2-\delta/2)_+}\right\rbrace\nonumber\\
&\leq& \!
\frac{2\lambda_2(f_\varepsilon,\kappa_0)\Gamma_2(m)}{n}\exp\left\lbrace-\frac{K'\xi^2}{2}(\pi
m)^{(1/2-\delta/2)_+}\right\rbrace\!\sum_{m'=1}^{m_n} \!\exp\left\lbrace-\frac{K'\xi^2}{2}(\pi m')^{(1/2-\delta/2)_+}\right\rbrace\nonumber \\
&&+
\sum_{m'=1}^{m_n}\frac{2\lambda_2(f_\varepsilon,\kappa_0)\Gamma_2(m')}{n}\exp\left\lbrace-\frac{K'\xi^2}{2}(\pi
m')^{(1/2-\delta/2)_+}\right\rbrace \label{I}.
\end{eqnarray}

\paragraph{\textbf{Case $0\leq\delta < 1/3$}} In that case, since
$\delta< (1/2-\delta/2)_+$,  the choice $\xi^2=1$ ensures that the
quantity  $\Gamma_2(m)\exp\{-(K'\xi^2/2)({m})^{(1/2-\delta/2)}\}$ is
bounded, and thus the first term in \eref{I} is bounded by $C/n.$
Since $1\leq m\leq m_n$ with $m_n$ satisfying \eref{mn},
$n^{-1}\sum_{m'=1}^{m_n}\Gamma_2(m')\exp\{-(K'/2)({m'})^{(1/2-\delta/2)}\}
$ is bounded by ${\tilde C}/n,$ and hence $\sum_{m'|m^*\geq
m_0}I(m^*)\leq Dn^{-1}.$ According to \eref{pmmp}, the result follows by
choosing $\mbox{pen}(m)=4ap(m,m')=24an^{-1}\Delta(m). $

\medskip

\paragraph{\textbf{Case $\delta= 1/3$}} According to
\eref{I}, we choose $\xi^2$  such that $2\mu\pi^{\delta}(m)^{\delta}
- (K'\xi^2/2){m}^{\delta}= -2\mu (\pi {m})^{\delta}$ that is
$\xi^2=(8\mu\pi^\delta\lambda_2(f_\varepsilon,\kappa_0))/(K_1\lambda_1(f_\varepsilon,\kappa_0')).$
Arguing as for the case $0\leq \delta<1/3$, this choice ensures that
$\sum_{m'|m^*\geq m_0}I(m^*) \leq Dn^{-1}$, and consequently \eref{but1}
holds. The result follows by taking $p(m,m')=2(1+2\xi^2)\Delta(m^*)n^{-1},$ and $\mbox{pen}(m)=8a(1+2\xi^2)\Delta(m)n^{-1}. $

\medskip

\paragraph{\textbf {Case $\delta> 1/3$}} In that case $\delta >
(1/2-\delta/2)_+$. Choose   $\xi^2(m)$ such that
$2\mu\pi^{\delta}(m)^{\delta} - (K'\xi^2/2){m}^{(1/2-\delta)_+}=
-2\mu\pi^{\delta}(m)^{\delta}$. Hence
$\xi^2(m)=(8\mu(\pi)^{\delta}\lambda_2(f_\varepsilon,\kappa_0)/(K_1\lambda_1(f_\varepsilon,\kappa_0'))
{(\pi m)}^{\delta-(1/2-\delta/2)_+}$. This choice ensures that
$\sum_{m'|m^*\geq m_0}I(m^*) \leq D/n$, so that \eref{but1} holds.
The result follows by choosing $p(m,m')=2(1+2\xi^2(m^*))\Delta(m^*)/
n,$ associated to $\mbox{pen}(m)=8a(1+2\xi^2(m))\Delta(m)/n.$

\medskip

\noindent\textbf{Proof of \eref{but2}.} Since $\max(m,\hat m)\leq
m_n$, according to \eref{lemnu},
\begin{eqnarray*}\sup_{t\in B_{m,\hat
m}(0,1)}\mathbb{E}\left(\nu_{n,X}(t)\right)^2 &\leq&\sup_{t \in
S_{m_n}, \|t\|=1} \mathbb{E}\left(\frac{1}{2\pi}\int
\nu_{n,X}(e^{ix\cdot})t^*(-x)dx\right)^2\\&\leq &
\frac{1}{2\pi}\int_{-\pi {m_n}}^{\pi
{m_n}}\mbox{Var}\left(\frac{1}{n}\sum_{k=1}^n
e^{ixX_k}\right)dx\\
&\leq& \frac{{m_n}}{n} +\frac{1}{\pi n}\int_{-\pi {m_n}}^{\pi {m_n}}
\sum_{k=2}^n\left\vert\mbox{Cov}\left(e^{ixX_1},e^{ixX_k}\right)\right\vert
dx
\end{eqnarray*}
and Theorem \ref{genepenmelss} is proved. \hfill $\Box$

\subsection{Proofs of Theorem \ref{genepenmelos} (1)} \label{secbeta}
We use the coupling argument recalled in Section \ref{couplage} to
build approximating variables for the $X_i$'s. For  $n=2p_nq_n+r_n$,
$0\leq r_n<q_n$, and $\ell =0, \cdots, p_n-1$, denote by
\begin{eqnarray*}
E_{\ell}=(X_{2 \ell q_n+1},...,X_{(2
  \ell+1)q_n}), && F_{\ell} =(X_{(2 \ell
  +1)q_n+1},...,X_{(2 \ell +2)q_n}), \\
E_{\ell}^{\star}=(X_{2 \ell q_n+1}^{\star},...,X_{(2
  \ell+1)q_n}^{\star}),&& F_{\ell}^{\star} =(X_{(2 \ell
  +1)q_n+1}^{\star},...,X_{(2 \ell +2)q_n}^{\star}).
\end{eqnarray*}
The variables $E_\ell^\star$ and $F_\ell^\star$ are such that
\begin{itemize}
\item[-] $E_{\ell}^{\star}$, $E_{\ell}$, $F_\ell^\star$ and $F_\ell$  are
identically distributed,
\item[-] $\mathbb{P}( E_{\ell}\neq E_{\ell}^{\star}) \leq
\beta_{\mathbf{X},\infty}(q_n)$ and
 $\mathbb{P}( F_{\ell}\neq
F_{\ell}^{\star})\leq \beta_{\mathbf{X},\infty}(q_n)$, \item[-] The
variables $(E_\ell^\star)_{0 \leq \ell \leq p_n-1}$ are i.i.d., and
so are the variables $(F_\ell^\star)_{0 \leq \ell \leq p_n-1}$.
\end{itemize}
 Without loss of
generality and for sake of simplicity we assume that $r_n=0$. For $\kappa_a$
defined in \eref{kappaa}, we start from
\begin{eqnarray*}
\|\tilde g-g\|^2 &\leq& \kappa_a^2\|g_m^{(n)} -g\|^2 +
2a\kappa_a\sup_{t\in B_{m,\hat m}(0,1)}(\nu_n^{(1)}(t))^2+
2a\kappa_a\sup_{t\in B_{m,\hat
m}(0,1)}(\nu_{n,X}(t))^2\\&&+\kappa_a({\rm pen}(m)- {\rm pen}(\hat
m))\\&\leq& \kappa_a^2\|g_m^{(n)} -g\|^2 + 2a\kappa_a\sup_{t\in
B_{m,\hat m}(0,1)}(\nu_n^{(1)}(t))^2+4a\kappa_a\sup_{t\in B_{m,\hat
m}(0,1)}(\nu_{n,X}^\star(t))^2\\&&+ 4a\kappa\sup_{t\in B_{m,\hat
m}(0,1)}(\nu_{n,X}(t)-\nu_{n,X}^\star(t))^2 +\kappa_a({\rm pen}(m)-
{\rm pen}(\hat m)),
\end{eqnarray*}
where $\nu_{n,X}^\star(t)$ is defined as $\nu_{n,X}(t)$ with
$X_i^\star$ instead of  $X_i$. Choose $p_1(m,m')$ and $p_2(m,m')$
such that
$$2ap_1(m,m')\leq [\mbox{pen}_1(m)+\mbox{pen}_1(m')]\mbox{ and
}4ap_2(m,m')\leq [\mbox{pen}_2(m)+\mbox{pen}_2(m')],$$ for
pen$(m)=\mbox{pen}_1(m)+\mbox{pen}_2(m).$
It follows that
\begin{eqnarray}
\|\tilde g-g\|^2&\leq&
\kappa_a^2\|g-g_m^{(n)}\|^2+2\kappa_a\mbox{pen}(m)+4a\kappa_aW_{n,X}^\star(m,\hat
m)+4a\kappa_a\sup_{t\in B_{m,\hat
m}(0,1)}(\nu_{n,X}(t)-\nu_{n,X}^\star(t))^2\nonumber\\&&+2a\kappa_a
W_n(m,\hat m)\nonumber\\&\leq&
\kappa_a^2\|g-g_m^{(n)}\|^2+2\kappa_a\mbox{pen}(m)+
4a\kappa_a\!\!\!\sum_{m'=1}^{m_n}W_{n,X}^\star(m,m') +2a\kappa_a
\sum_{m'=1}^{m_n}W_n(m,m')\\\nonumber&&+4a\kappa_a\sup_{t\in
B_{m,\hat
m}(0,1)}(\nu_{n,X}(t)-\nu_{n,X}^\star(t))^2,\label{departos}
\end{eqnarray}
where
\begin{eqnarray}\label{Wgos} W_n(m,m')&:=& \Big [\sup_{t\in B_{m, m'}(0,1)}
 |\nu_n^{(1)}(t)|^2-p_1(m, m')\Big ]_+,\\
W_{n,X}^\star(m,m')&:=&\Big [\sup_{t\in B_{m, m'}(0,1)}
 |\nu_{n,X}^\star(t)|^2-p_2(m, m')\Big ]_+.
\end{eqnarray}
The main parts of the proof lies in the three following points~:

\textbf{1)} Study of $W_n(m,m')$. More precisely, we have to find
$p_1(m,m')$ such that for a constant $A_2$,
\begin{equation}\label{but1beta}
\sum_{m'=1}^{m_n} \mathbb{E}(W_n(m,m'))\leq \frac{A_2}{n}.
\end{equation}

\textbf{2)} Study of $W_{n,X}^\star(m,m')$. More precisely, we
have to find $p_2(m,m')$ such that for a constant $A_3$,
\begin{equation}\label{but3beta}
\sum_{m'=1}^{m_n} \mathbb{E}(W_{n,X}^\star(m,m'))\leq \frac{A_3}{n}.
\end{equation}

\textbf{3)} Study of $\sup_{t\in B_{m,\hat
m}(0,1)}(\nu_{n,X}(t)-\nu_{n,X}^\star(t))^2$ and more precisely we
have to prove that
\begin{eqnarray}
\label{but2beta} \mathbb{E}\Big[\sup_{t\in B_{m,\hat
m}(0,1)}(\nu_{n,X}^\star(t)-\nu_{n,X}(t))^2\Big]\leq
4\beta_{\mathbf{X},\infty}(q_n){m_n}\leq \frac{A_4}{n}.
\end{eqnarray}

\noindent \textbf{Proof of \eref{but1beta}} The proof of
\eref{but1beta} for ordinary smooth errors ($\delta=0$ in
\eref{condfeps}) is the same as the proof of \eref{but1} by taking
$p_1(m,m')=p(m,m')$, with $p(m,m')$ as in the proof of \eref{but1}
and $\xi^2=1$. Hence we choose pen$_1(m)=24an^{-1}\Delta(m)$.\\

\noindent \textbf{Proof of \eref{but3beta}} We proceed as in the
independent case by applying Lemma \ref{Concent}. Set
$m^*=\max(m,m')$. The process $W_{n,X}^\star(m,m')$ must be split
into two terms $(W_{n,1,X}^\star(m,m')+W_{n,2,X}^\star(m,m'))/2$
involving respectively the odd and even blocks,  which are of the
same type. More precisely $W_{n,k,X}^{\star}(m,m')$ is defined,
for $k=1,2$, by
\begin{eqnarray*}
W_{n,k,X}^\star(m,m')= \Big[\sup_{t\in B_{m, m'}(0,1)}
    \Big|\frac{1}{p_nq_n}\sum_{\ell=1}^{p_n}\sum_{i=1}^{q_n}\Big(t(X^\star_{(2\ell+k-1)
      q_n+i})-\langle t, g\rangle \Big)\Big|^2-p_{2,k}(m,
    m')\Big]_+
\end{eqnarray*}
We only study $W_{n,1,X}^\star(m,m')$ and conclude for
$W_{n,2,X}^\star(m,m')$ by using analogous arguments.  The study of
$W_{n,1,X}^\star(m,m')$ consists in applying  Lemma \ref{Concent}
to $\nu_{n,1,X}^{\star}(t)$ defined by
\begin{eqnarray*}
 \nu_{n,1,X}^{\star}(t)
=\frac{1}{p_n}\sum_{\ell
  =1}^{p_n}\nu_{q_n,\ell,X}^\star(t) \mbox{ with }\nu_{q_n,\ell,X}^{\star}(t)=\frac{1}{q_n}\sum_{j=1}^{q_n}t(X^{\star}_{2\ell
q_n
  +j})-\langle t, g\rangle,
\end{eqnarray*}
considered as the sum of the $p_n$ independent random variables
$\nu_{q_n,\ell,X}^{\star}(t)$.
Denote by $M_1^\star(m^*)$, $H^\star(m^*)$ and $v^\star(m^*)$ 
quantities such that
\begin{eqnarray*}
\sup_{t \in
  B_{m,m'}(0,1)}\parallel\nu_{q_n,\ell,X}^{\star}(t)\parallel_\infty & \leq &
M_1^\star(m^*), \quad  \mathbb{E}\Big (\sup_{t \in
  B_{m,m'}(0,1)}\vert\nu_{n,1,X}^{\star}(t)\vert \Big )\leq H^\star(m^*)\\
     &\text{and}&   \sup_{t
  \in B_{m,m'}(0,1)}\mbox{Var}( \nu_{q_n,\ell,X}^{\star}(t) ) \leq
v^\star(m^*). \end{eqnarray*}
Lemma \ref{l2} leads to the choices $M_1^\star(m^*)=\sqrt{m^*}$,
\begin{eqnarray*}
(H^{\star}(m^*))^2=\frac{\Big(1+4
\sum_{k=1}^n\beta_{\mathbf{X},1}(k)\Big){m^*}}{n}, \quad
 \text{ and }v^\star(m^*)=\frac{8\Big(
\sum_{k=0}^{q_n}(k+1)\beta_{\mathbf{X},1}(k)\|g\|_{\infty} m^*\Big)^{1/2}}{q_n}.\end{eqnarray*}
Take $\xi^2(m^*)=1/2$. We use that for  $m^*\geq
m_0$, $$2(1+2\xi^2(m^*))(H^\star(m^*))^2=4(H^\star(m^*))^2\leq \Delta(m^*)/(4n).$$
Then we take
$p_{2,1}(m,m')=\Delta(m)/(4n),$ 
and get that
\begin{eqnarray*}
\sum_{m'=1}^{m_n}\mathbb{E}(W_{n,1,X}^\star(m,m'))&=&\sum_{m'|m^*\leq m_0}
\mathbb{E}(W_{n,1,X}^\star(m,m'))+\sum_{m'|m^*>m_0}\mathbb{E}(W_{n,1,X}^\star(m,m'))\\
&\leq&
\sum_{m'|m^*\leq m_0}\mathbb{E}\Big[\sup_{t\in
    B_{m,m'}(0,1)}|\nu_{n,1,X}^\star(t)|^2-4(H^\star(m^*))^2\Big]_+\\&&
+\sum_{m'|m^*\leq m_0}\vert
p_{21}(m,m')-4 (H^\star(m^*))^2\vert\\&&+
\sum_{m'|m^*>m_0}\mathbb{E}\Big[\sup_{t\in
    B_{m,m'}(0,1)}|\nu_{n,1,X}^\star(t)|^2-4(H^\star(m^*))^2\Big]_+.
\end{eqnarray*}
It follows that
\begin{eqnarray*}
\sum_{m'=1}^{m_n}\mathbb{E}(W_{n,1,X}^\star(m,m'))&\leq& 2\sum_{m'=1}^{m_n}\mathbb{E}\Big[\sup_{t\in
    B_{m,m'}(0,1)}|\nu_{n,1,X}^\star(t)|^2-4(H^\star(m^*))^2\Big]_+\\&&
+\sum_{m'|m^*\leq m_0}\vert
p_{2,1}(m,m')-4 (H^\star(m^*))^2\vert\\&\leq&
2\sum_{m'=1}^{m_n}\mathbb{E}\Big[\sup_{t\in
    B_{m,m'}(0,1)}|\nu_{n,1,X}^\star(t)|^2-4(H^\star(m^*))^2\Big]_+
+\frac{C(m_0)}{n}.
\end{eqnarray*}
We apply Lemma \ref{Concent} to $\mathbb{E}\Big[\sup_{t\in
    B_{m,m'}(0,1)}|\nu_{n,1,X}^\star(t)|^2-4(H^\star(m^*))^2\Big]_+$
and obtain
$$\sum_{m'=1}^{m_n}\mathbb{E}\Big[\sup_{t\in
    B_{m,m'}(0,1)}|\nu_{n,1,X}^\star(t)|^2-4(H^\star(m^*))^2\Big]_+ \leq K
\sum_{m'=1}^{m_n}[I^\star(m^*)+II^\star(m^*)],$$ with $I^\star(m^*)$
and $II^\star(m^*)$ defined by
\begin{eqnarray*}
I^{\star}(m^*)=\frac{m^*}{n} \exp\Big \{-K_2\sqrt{m^*}\Big
\} \mbox{ and } II^{\star}(m^*)= \frac{ q_n^2
{m^*}}{n^2}\exp\Big\{-\frac{\sqrt{2}K_1 \xi
C(\xi)}{7}\frac{\sqrt{n}}{q_n} \Big\},
\end{eqnarray*}
where $K_2=(K_1/32)(1+4\sum_{k=1}^n \beta_{{\bf X},1}(k))/\sqrt{\|g\|_{\infty} \sum_{k=0}^{q_n}(k+1)\beta_{{\bf X},1}(k)}$.\\
With our choice of $\xi^2(m)$, if we take $q_n=[n^c]$, for $c$
in $]0,1/2[$, then
$$\sum_{m'}I(m^*)\leq \frac{C}{n},\mbox{ and }
\sum_{m'=1}^{m_n}II^\star(m^*)\leq \frac{C}{n}.$$
Finally $$\sum_{m'=1}^{m_n}\mathbb{E}\Big[\sup_{t\in
    B_{m,m'}(0,1)}|\nu_{n,1,X}^\star(t)|^2-4(H^\star(m^*))^2\Big]_+\leq\frac{C}{n}
$$ and
$$\sum_{m'=1}^{m_n}\mathbb{E} [W_{n,X}^\star(m,m')]\leq
2\sum_{m'=1}^{m_n}\mathbb{E} [W_{n,1,X}^\star(m,m')+
W_{n,2,X}^\star(m,m')] \leq \frac{C}{n}.$$ The result follows for
choosing $p_2(m,m')=2p_{2,1}(m,m')+2p_{2,2}(m,m')=\Delta(m)/n, $
\mbox{ and }\mbox{pen}$(m)=25a\Delta(m)/n.$

\medskip

\noindent \textbf{Proof of \eref{but2beta}.} A rough bound is
obtained by writing that
\begin{eqnarray*}
\sup_{t\in B_{m,\hat
m}(0,1)}\vert\nu_{n,X}^\star(t)-\nu_{n,X}(t)\vert^2&=&\sup_{t\in
S^{(n)}_{\max(m,\hat m)}, \parallel t\parallel\leq
1}\vert\nu_{n,X}^\star(t)-\nu_{n,X}(t)\vert^2 \\ &\leq&
 \sup_{t\in S_{m_n}, \parallel t\parallel\leq
1}\vert\nu_{n,X}^\star(t)-\nu_{n,X}(t)\vert^2.
\end{eqnarray*}
According to \eref{lemnu},
\begin{eqnarray*}
\nu_{n,X}^\star (t)-\nu_{n,X}(t) \!\!\!&=&\!\!\!\frac{1}{2\pi}\int
[\nu_{n,X}^\star(e^{ix\cdot})- \nu_{n,X}(e^{ix\cdot})]t^*(-x)dx.
\end{eqnarray*}
Since $\vert
\nu_{n,X}(e^{ix\cdot})-\nu_{n,X}^\star(e^{ix\cdot})\vert \leq 2$, we
have
\begin{eqnarray*}
\sup_{t\in B_{m,\hat
m}(0,1)}\vert\nu_{n,X}^\star(t)-\nu_{n,X}(t)\vert^2
&\leq &
 \sup_{t\in S_{m_n}, \parallel t\parallel\leq 1}
\frac{1}{4\pi^2}\left|\int
[\nu_{n,X}^\star(e^{ix\cdot})- \nu_{n,X}(e^{ix\cdot})]t^*(-x)dx\right|^2\\
&\leq& \frac{1}{2\pi}\int_{-\pi {m_n}}^{\pi {m_n}}
\vert \nu_{n,X}^\star(e^{ix\cdot})- \nu_{n,X}(e^{ix\cdot})\vert^2dx\\
&\leq & \frac{1}{\pi}\int_{-\pi {m_n}}^{\pi {m_n}} \vert
\nu_{n,X}^\star(e^{ix\cdot})- \nu_{n,X}(e^{ix\cdot})\vert dx.
\end{eqnarray*}
According to the properties of the coupling,
\begin{eqnarray*}
\mathbb{E}\Big[\sup_{t\in B_{m,\hat
m}(0,1)}\vert\nu_{n,X}^\star(t)-\nu_{n,X}(t)\vert^2 \Big] &\leq &
\frac{1}{\pi}  \int_{-\pi {m_n}}^{\pi {m_n}} \mathbb{E}\vert
\nu_{n,X}^\star(e^{ix\cdot})- \nu_{n,X}(e^{ix\cdot})\vert dx\leq
4\beta_{\mathbf{X},\infty}(q_n) {m_n}.
\end{eqnarray*}
For ordinary smooth errors, according to \eref{mn}, ${m_n}\leq
n^{1/(2\gamma+1)}$. It follows that if we choose $q_n$ such that
$\beta_{\mathbf{X},\infty}(q_n)=O( n^{-(2\gamma+2)/(2\gamma+1)})$,
then $\beta_{\mathbf{X},\infty}(q_n){m_n}=O(n^{-1}) $.
For $q_n=[n^c]$ and $\beta_{\mathbf{X},\infty}(n)=
O(n^{-1-\theta})$, we obtain the condition
$n^{-c(1+\theta)}=O( n^{-(2\gamma+2)/(2\gamma+1)})$. If
$\theta>(2\gamma+3)/(2\gamma+1)$, one can find $c<1/2$ such that
this condition is satisfied.

\subsection{Proofs of Theorem \ref{genepenmelos} (2)}
We proceed as in the $\beta$-mixing case, by using the coupling
argument given in Section \ref{couplage}. The variables $E_\ell,
E_\ell^\star, F_\ell, F_\ell^\star $ are build as in Section
\ref{secbeta} and are such that
\begin{itemize}
\item[-] $E_{\ell}^{\star}$, $E_{\ell}$, $F_\ell^\star$ and $F_\ell$  are
identically distributed,
\item[-] $\displaystyle \sum_{i=1}^{q_n} {\mathbb E} (|X_{2 \ell q_n+i}-X_{2 \ell
  q_n+i}^\star|) \leq q_n
\tau_{\mathbf{X},\infty}(q_n)$ and
$\displaystyle \sum_{i=1}^{q_n}
{\mathbb E} (|X_{(2 \ell
  +1)q_n+i}-X_{(2 \ell
  +1)q_n+i}^\star|) \leq q_n
\tau_{\mathbf{X},\infty}(q_n)$, \item[-] The variables
$(E_\ell^\star)_{0 \leq \ell \leq p_n-1}$ are i.i.d., and so are the
variables $(F_\ell^\star)_{0 \leq \ell \leq p_n-1}$.
\end{itemize}
Without
loss of generality and for sake of simplicity we assume that
$r_n=0$.
As for the proof of Theorem \ref{genepenmelos} under \textbf{2)}, we
start from \eref{departos}. Hence we have to~:

\medskip

\textbf{1)} Study of $W_n(m,m')$, and more precisely in finding
$p_1(m,m')$ such that for a constant $K_2$,
\begin{equation}\label{but1tau}
\sum_{m'=1}^{m_n} \mathbb{E}(W_n(m,m'))\leq \frac{K_2}{n}.
\end{equation}

\textbf{2)} Study of $W_{n,X}^\star(m,m')$, and more precisely in
finding $p_2(m,m')$ such that for a constant $K_3$,
\begin{equation}\label{but3tau}
\sum_{m'=1}^{m_n} \mathbb{E}(W_{n,X}^\star(m,m'))\leq \frac{K_3}{n}.
\end{equation}

\textbf{3)} Study of $\sup_{t\in B_{m,\hat
m}(0,1)}(\nu_{n,X}(t)-\nu_{n,X}^\star(t))^2$ and more precisely in proving that
\begin{eqnarray}
\label{but2tau} \mathbb{E}\big[\sup_{t\in B_{m,\hat
m}(0,1)}(\nu_{n,X}^\star(t)-\nu_{n,X}(t))^2\big]\leq \pi
\tau_{\mathbf{X},\infty}(q_n){m_n}^2\leq \frac{K_4}{n}.
\end{eqnarray}

\medskip

\noindent \textbf{Proof of \eref{but1tau}} The proof of
\eref{but1tau} for ordinary smooth errors is the same as the proof
of \eref{but1}.

\medskip

\noindent \textbf{Proof of \eref{but3tau}} As for the proof
\eref{but3beta} we apply Lemma \ref{Concent} with
\begin{eqnarray*}
(H^{\star}(m^*))^2=
\frac{\left(m^*+\pi\sum_{k=1}^{n-1}\tau_{\mathbf{X},1}(k)(m^*)^2\right)}{n},
\quad M_1^\star(m^*)=m^*,\end{eqnarray*}
\begin{eqnarray*}
\mbox{ and }
v^\star(m^*)=\frac{\left(m^*+\pi\sum_{k=1}^{n-1}\tau_{\mathbf{X},1}(k)(m^*)^2\right)}{q_n}.
\end{eqnarray*}
We take $\xi^2=\xi^2(m)=(3/K_1+1)\ln(m)$.
 In the same way as for the proof Theorem \ref{genepenmelos}(1), we use that for  $m^*\geq m_0$,
$$2(1+2\xi^2(m^*))(H^\star(m^*))^2\leq \Delta(m^*)/(4n).$$ Then we take
$p_{21}(m,m')=\Delta(m^*)(4n)^{-1}$ and  get that
\begin{eqnarray*}
\sum_{m'=1}^{m_n}\mathbb{E}(W_{n,1,X}^\star(m,m'))\leq
2\sum_{m'=1}^{m_n}\mathbb{E}\Big[\sup_{t\in
    B_{m,m'}(0,1)}|\nu_{n,1,X}^\star(t)|^2-2(1+2\xi^2(m^*))(H^\star(m^*))^2\Big]_+
+\frac{C(m_0)}{n}.
\end{eqnarray*}
We now apply Lemma \ref{Concent} to $\mathbb{E}\Big[\sup_{t\in
    B_{m,m'}(0,1)}|\nu_{n,1,X}^\star(t)|^2-2(1+2\xi^2(m^*))(H^\star(m^*))^2\Big]_+$ and obtain
$$\sum_{m'=1}^{m_n}\mathbb{E}\Big[\sup_{t\in
    B_{m,m'}(0,1)}|\nu_{n,1,X}^\star(t)|^2-2(1+2\xi^2(m^*))(H^\star(m^*))^2\Big]_+ \leq K
\sum_{m'=1}^{m_n}[I^\star(m^*)+II^\star(m^*)],$$
 with $I^\star(m^*)$
and $II^\star(m^*)$ now defined by
\begin{eqnarray*}
&&I^{\star}(m^*)=\frac{{m^*}^2}{n}
\exp\{-K_1\xi^2(m^*)\}\\
\mbox{ and } &&II^{\star}(m^*)= \frac{ q_n^2
{m^*}^2}{n^2}\exp\left\lbrace -\frac {\sqrt{2}K_1 \xi
C(\xi)\Big(1+\pi\sum_{k=1}^n\tau_{\mathbf{X},1}(k)\Big)}{7}\frac{\sqrt{n}}{q_n}\right\rbrace.
\end{eqnarray*}
With this $\xi^2(m)$, if we take  $q_n=[n^c]$, with $c$ in $]0,1/2[$
then $$\sum_{m'}I(m^*)\leq \frac{C}{n} \quad \mbox{ and } \quad
\sum_{m'=1}^{m_n}II(m^*)
\leq \frac{C}{n}.$$ Finally $\sum_{m'=1}^{m_n}\esp
[W_n^\star(m,m')]\leq 2\sum_{m'=1}^{m_n}\mathbb{E}
[W_{n,1,X}^\star(m,m')+ W_{n,2,X}^\star(m,m')] \leq C n^{-1}.$ The
result follows by choosing
$p_2(m,m')=2p_{21}(m,m')+2p_{22}(m,m')=\Delta(m) n^{-1}$, and
$\mbox{pen}(m)=25a\Delta(m) n^{-1}.$

\medskip

\noindent \textbf{Proof of \eref{but2tau}} The proof of
\eref{but2tau} is  similar to the proof of  \eref{but2}. Since
$|e^{-ix t}-e^{-ix s}| \leq |x| |t-s|$, one has
$$
 \sum_{i=1}^{q_n} {\mathbb E} (|e^{-i X_{2 \ell q_n+i}}-e^{-i X_{2 \ell
  q_n+i}^\star}|) \leq q_n |x|
\tau_{\mathbf{X},\infty}(q_n)
$$
It follows that
\begin{eqnarray*}
\mathbb{E}\Big [\sup_{t\in B_{m,\hat
m}(0,1)}\vert\nu_{n,X}^\star(t)-\nu_{n,X}(t)\vert^2 \Big ] &\leq &
\frac{1}{\pi}\int_{-\pi m_n}^{\pi m_n} \mathbb{E}\vert
\nu_{n,X}^\star(e^{ix\cdot})- \nu_{n,X}(e^{ix\cdot})\vert dx\leq
 \pi \tau_{\mathbf{X},\infty}(q_n) {m_n}^2.
\end{eqnarray*}
For ordinary smooth errors, according to \eref{mn}, ${m_n}^2\leq
n^{2/(2\gamma+1)}$. It follows that if we choose $q_n$ such that
$\tau_{\mathbf{X},\infty}(q_n)=O( n^{-(2\gamma+3)/(2\gamma+1)})$,
then $\tau_{\mathbf{X},\infty}(q_n){m_n^2}=O(n^{-1}) $.
For $q_n=[n^c]$ and $\tau_{\mathbf{X},\infty}(n)= O(n^{-1-\theta})$,
we obtain the condition
$n^{-c(1+\theta)}=O( n^{-(2\gamma+3)/(2\gamma+1)})$. If
$\theta>(2\gamma+5)/(2\gamma+1)$, one can find $c<1/2$ such that
this condition is satisfied.
\hfill $\Box$

\subsection{Proof of Corollary \ref{densitycor}}
The result follows from the proof of Theorem \ref{genepenmelos} (1),
where only the process $\nu_{n,X}$ appears.$\Box$

\section{Technical  lemmas}
\setcounter{equation}{0}
\setcounter{lem}{0}
\setcounter{thm}{0}

\begin{lem}
\label{nucond}If we denote by $\nu_{n,X}(t)$ the quantity defined by
\eref{nu}, then
\begin{eqnarray*}
n^{-1}\sum_{k=1}^n\mathbb{E}(u_t^*(Z_k)|\sigma(X_i,\,i\geq
0))-<t,g>=\nu_{n,X}(t).
\end{eqnarray*}
\end{lem}
The proof of Lemma \ref{nucond}, rather straightforward, is omitted.
\begin{lem}
\label{details} Let $\nu_{n,Z}(u_t^*)$ be defined by \eref{nu},
$\Delta(m)$ being defined in \eref{Delta1}. Then
\begin{equation*}
\sum_{j \in \mathbb{Z}} \left\vert u^*_{\varphi_{m,j}}(z)\right\vert
^2=(2\pi)^{-1} m\int \left\vert
\frac{\varphi^*(x)}{f_\varepsilon^*(xm)}\right\vert^2dx= \Delta(m).
\end{equation*}
\end{lem}


\begin{lem}
\label{l3} Let $\nu_{n,Z}(u_t^*)$, $\Delta(m)$ and $\Delta_2(m, h)$
be defined in \eref{nu}, \eref{Delta1} and in \eref{Delta2}. Then
\begin{eqnarray*}
\sup_{ t \in B_{m,m'}(0,1) }\parallel u_t^*\parallel_\infty \leq
\sqrt{\Delta(m^*)}&&  \mathbb{E}[\sup_{t\in
B_{m,m'}(0,1)}|\nu_{n,Z}(u_t^*)|]\leq \sqrt{\Delta(m^*)/n},\\
{\rm and } \sup_{t\in B_{m,m'}(0,1)} {\rm Var}(u^*_t(Z_1))&\leq &
\sqrt{\Delta_2(m^*,h)}/(2\pi).
\end{eqnarray*}
\end{lem}

\noindent We refer to Comte {\it et al.} \citeyear{CRT1} for the
proofs of Lemmas \ref{details} and \ref{l3}.

\medskip

\begin{lem}
\label{norminf}
$\parallel\sum_{j\in \mathbb{Z}} |\varphi_{m,j}|^2\parallel_\infty \leq
m.$
\end{lem}
\noindent \textbf{Proof of Lemma \ref{norminf}}
Write
\begin{eqnarray*}
\sum_{j \in \mathbb{Z}} \left\vert \varphi_{m,j}(x)\right\vert ^2
&=&
\frac{1}{(2\pi)^2}\sum_{j \in \mathbb{Z}} \left\vert \int e^{-iux} \varphi_{m,j}^*(u)
  du\right\vert ^2
=
\frac{m}{(2\pi)^2}\sum_{j \in \mathbb{Z}} \left\vert \int
e^{-ixum}e^{iju}\varphi^*(u)du\right\vert^2.
\end{eqnarray*}
We conclude by applying Parseval's Formula
which gives that
\begin{equation*}
\sum_{j \in \mathbb{Z}} \left\vert
\varphi_{m,j}(x)\right\vert ^2=(2\pi)^{-1} m\int
\left\vert
\varphi^*(u)\right\vert^2du= m.
\end{equation*}

\begin{lem}
\label{l2} For $B_{m, m'}(0,1)=\{t\in S_{m\vee m'} \;/\;
\|t\|_2=1\}$,  we have, for $m^*=m\vee m'$,
\begin{eqnarray*}
\sup_{ t \in B_{m,m'}(0,1) }\parallel t\parallel_\infty  \leq
\sqrt{m^*},~\hspace{0.5cm}~\mathbb{E}[\sup_{t\in
B_{m,m'}(0,1)}|\nu_{n,1,X}^{\star}(t)|]\leq
\sqrt{\frac{(1+4\sum_{k=1}^n \beta_{{\bf X},1}(k)) m^*}n}\\\mbox{
and } \sup_{t\in B_{m,m'}(0,1)} {\rm Var}(\nu_{q_n, \ell,
X}^{\star}(t))\leq \frac{\left[2 \|g\|_{\infty}(1+ 32\sum_{k=1}^n
(1+k)\beta_{{\bf X}, 1} (k)) \right]^{1/2} \sqrt{m^*}}{q_n}.
\end{eqnarray*}
\end{lem}

\noindent \textbf{Proof of Lemma \ref{l2}} For $t$ in $B_{m,
m'}(0,1)$, with $m^*=m\vee m'$, one has $t=\sum_{j \in
  \mathbb{Z}}b_{m^*,j}\varphi_{m^*,j}$. Applying
Cauchy-Schwarz Inequality and Lemma \ref{norminf} we obtain
\begin{eqnarray*}
\sup_{ t \in B_{m,m'}(0,1) }\parallel t\parallel_\infty \leq \Big\|
\sum_{j\in \mathbb{Z}} |\varphi_{m^*,j}|^2 \Big\|_\infty^{1/2}
  \leq \sqrt{m^*}.
\end{eqnarray*}
Now,  using again Cauchy-Schwarz Inequality
\begin{eqnarray*}
\mathbb{E}\left[\sup_{t\in
B_{m,m'}(0,1)}|\nu_{n,1,X}^{\star}(t)|\right]&\leq&
\esp\left[\sqrt{\sum_{j\in \mathbb{Z}} (\nu_{n,1 X}^{\star}
(\varphi_{m^*,j}))^2}\right] \leq\sqrt{\sum_{j\in \mathbb{Z}
  }\mbox{Var}(\nu_{n,1,X}^{\star}(\varphi_{m^*,j}))}.
  \end{eqnarray*}
By analogy with (\ref{varineq}), we write
$$
\mathbb{E}\left(\sum_{j\in
\mathbb{Z}}\left(\nu_{n,1,X}^{\star}(\varphi_{m,j})\right)^2\right)=\frac{1}{4\pi^2}\sum_{j\in
\mathbb{Z}}\mathbb{E}\Big(\int
\varphi_{m,j}^*(-x) \nu_{n,1,X}^{\star}(e^{ix\cdot})dx\Big)^2=\frac{1}{2\pi}\int_{-\pi
m}^{\pi m}\mathbb{E}\vert
\nu_{n,1,X}^{\star}(e^{ix\cdot})\vert^2dx.
$$
This yields
$$\mathbb{E}\left[\sup_{t\in B_{m,m'}(0,1)}|\nu_{n,1,X}^{\star}(t)|\right]\leq \frac{\left(1+4\sum_{k=1}^n
\beta_{{\bf X},1}(k)\right)m^*}n.$$ Finally, we apply
Viennet's~\citeyear{Viennet} variance inequality (see  Theorem 2.1
p. 472 and  Lemma 4.2 p. 481). Hence there exist some measurable
functions $b_k$,  such that   $0\leq b_k\leq 1$ and ${\mathbb
E}\left[\left(\sum_{k=1}^n b_k(X_1)\right)^2\right]\leq \sum_{k\geq
1}(1+k)\beta_{{\bf X}, 1}(k)$, for which
$$
\sup_{ t\in B_{m,m'}(0,1)} \mbox{Var}(\nu_{q_n,\ell,X}(t))\leq
\sup_{ t\in
  B_{m,m'}(0,1)}\frac 1{q_n} \int \left(1+4\sum_{k=1}^{q_n}b_k\right) t^2(x) g(x)dx \, .
  $$
Consequently
\begin{eqnarray*}
 \sup_{ t\in B_{m,m'}(0,1)} \mbox{Var}(\nu_{q_n,\ell,X}(t)) &\leq & \sup_{ t\in B_{m,m'}(0,1)}\frac 1{q_n}
\|t\|_{\infty} \|g\|_{\infty}^{1/2} \left[\int
\left(1+4\sum_{k=1}^{q_n}b_k\right)^2  g(x)dx\right]^{1/2}\\ &\leq &
\sqrt{2\|g\|_{\infty}(1+32\sum_{k=1}^{q_n}(1+k) \beta_{{\bf
X},1}(k))} \frac{\sqrt{m^*}}{q_n}.
\end{eqnarray*}


\noindent \textbf{Proof of Lemma \ref{Concent}~:} Starting from the
concentration inequality given in Klein and Rio \citeyear{KleinRio}
and arguing as in Birg\'e and Massart \citeyear{BirgeMassart98} (see
the  proof of their Corollary 2 page 354) we obtain the upper bound
\begin{eqnarray}
\label{TalKR}
\mathbb{P}\left( \sup_{g\in {\mathcal
G}}|\nu_n(g)|\geq (1+\eta)H+\lambda  \right)\leq
2\exp\left[-K_1n\left(\frac{\lambda^2}{v}\wedge\frac{2\lambda(\eta \wedge
1)}{7M_1}\right)   \right],
\end{eqnarray}
 where $K_1=1/6.$
By taking $\eta=(\sqrt{1+\epsilon}-1)\wedge 1=C(\epsilon)\leq 1$
we get
\begin{eqnarray*}
\mathbb{E}[\sup_{g\in {\mathcal
G}}|\nu_n(g)|^2-2(1+2\epsilon)H^2]_+&\!\!\!\!\!\!\leq &\!\!\!\int_0^{+\infty}
\!\!\!\!\!\!\mathbb{P}\left( \sup_{g\in {\mathcal
G}}|\nu_n(g)|^2\geq 2(1+2\epsilon)H^2+\tau\right)d\tau\\
\!\!\!\!\!\!\!\!\!&\leq &\!\!\!\int_0^{+\infty}
\!\!\!\!\!\!\mathbb{P}\left( \sup_{g\in {\mathcal
G}}|\nu_n(g)|\geq \sqrt{2(1+\epsilon)H^2+2(\epsilon H^2+\tau/2)}\right)d\tau\\
\!\!\!\!\!\!\!\!\!&\leq &\!\!\!2\int_0^{+\infty}
\!\!\!\!\!\!\mathbb{P}\left( \sup_{g\in {\mathcal
G}}|\nu_n(g)|\geq \sqrt{(1+\epsilon)}H+\sqrt{\epsilon H^2+\tau/2}\right)d\tau\\
\!\!\!\!\!\!\!\!\!&\leq
&\!\!\!4\left(\int_0^{+\infty}e^{-\frac{K_1n}{v}(\epsilon
H^2+\tau/2)} d\tau + \int_0^{+\infty}
e^{-\frac{2K_1nC(\epsilon)}{7M_1\sqrt{2}}(\sqrt{\epsilon}
H+\sqrt{\tau/2})}d\tau \right)\\
\!\!\!\!\!\!\!\!\!&\leq
&\!\!\!4e^{-K_1\epsilon\frac{nH^2}{v}}\int_0^{+\infty}\!\!\!\!e^{-\frac{K_1n}{2v}\tau}d\tau
+ 4e^{-\frac{\sqrt{2}K_1C(\epsilon
)\sqrt{\epsilon}}{7}\frac{nH}{M_1}} \int_0^{+\infty}\!\!\!\!\
e^{-\frac{ K_1C(\epsilon)
    n\sqrt{\tau}}{7 M_1}}d\tau.
\end{eqnarray*}
Using that for any positive constant $C$,
$\int_0^{+\infty}e^{-Cx}dx=1/C$ and
$\int_0^{+\infty}e^{-C\sqrt{x}}dx=2/C^2$, we get that
\begin{eqnarray*}
\mathbb{E}[\sup_{g\in {\mathcal
G}}|\nu_n(g)|^2-2(1+2\epsilon)H^2]_+&\!\!\!\!\leq &\!\!\!
\frac{8}{K_1}\left(\frac{v}{n}e^{-K_1\epsilon\frac{ nH^2}{v}}
    +\frac{49M_1^2}{K_1^2n^2C^2(\epsilon)}
e^{-\frac{\sqrt{2} K_1 C(\epsilon)\sqrt{\epsilon}
}{7}\frac{nH}{M_1}}  \right).\qquad\hfill \Box
\end{eqnarray*}

\bibliographystyle{plain}
\bibliography{Bibliodeconvmix}

\end{document}